\documentclass[10pt,a4paper]{amsart}
\usepackage [english]{babel}
\usepackage{hyperref}
\usepackage[usenames]{color}
\usepackage{amsthm}
\usepackage{longtable}
\usepackage{amsfonts, amssymb, amsmath, amscd, latexsym}
\usepackage{mathrsfs,epsfig}
\usepackage[latin1]{inputenc}
\usepackage{caption}

\newtheorem{thm}{Theorem}[section]

\newtheorem{theorem}{Theorem}

\newtheorem{lemma}[thm]{Lemma}

\newtheorem{prop}[thm]{Proposition}

\theoremstyle{remark}
\newtheorem{remark}[thm]{Remark}

\newtheorem{defin}[thm]{Definition}
\numberwithin{equation}{section}
\def\ds{\displaystyle}

\newcommand{\RR}{{\mathbb R}}
\newcommand{\lir}{{\lim_{r \to +\infty}}}
\newcommand{\liro}{{\lim_{r \to 0}}}

\def\la{\lambda}
\def\R{\boldsymbol{R}}
\def\Q{\boldsymbol{Q}}
\def\eu{\textrm{e}}
\def\ep{\varepsilon}

\def\PP{\mathcal{P}}
\def\Id{\mathbb{I}}

\newcommand{\vt}[1]{{\left\vert\kern-0.25ex\left\vert\kern-0.25ex\left\vert #1
    \right\vert\kern-0.25ex\right\vert\kern-0.25ex\right\vert}}

\newcommand\bs[1]{{\boldsymbol{#1}}}

\newcommand\ov[1]{\overline{#1}}
\newcommand\und[1]{\underline{#1}}

\DeclareMathSizes{5}{5}{3}{2} 
\date{\today}
\begin{document}

\author{Luca Bisconti} \author{ Matteo Franca} \address[L.\
Bisconti]{Dipartimento di Matematica e Informatica ``U. Dini",
  Universit\`a degli  Studi di Firenze, Via S.\ Marta 3, I-50139
  Firenze, Italy.   Partially supported by G.N.A.M.P.A. - INdAM
  (Italy).  }
  \email{luca.bisconti@unifi.it}
  \address[M.\ Franca]{Dipartimento di Ingegneria Industriale e Scienze Matematiche,
  Universit\`a Politecnica delle Marche, Via Brecce Bianche, I-60131
  Ancona, Italy. Partially supported by
  G.N.A.M.P.A. - INdAM (Italy) }
   \email{franca@dipmat.univpm.it}

\subjclass[2010]{Primary 35k58, 35k91, 34e05; Secondary   35b08, 35b35} 

\title{On a nonlinear parabolic problem: Stability properties of Ground States}

\begin{abstract}
  We consider the Cauchy-problem for the following parabolic equation:
  \begin{equation*}
    \displaystyle u_t = \Delta u+ f(u,|x|),
  \end{equation*}
  where $x \in \RR^n$, $n >2$, and $f=f(u,|x|)$ is either critical or
  supercritical with respect to the Joseph-Lundgren exponent.  Using a
  new unifying approach we extend to a larger class of nonlinear
  potentials $f$, some known results concerning stability and weak
  asymptotic stability of positive Ground States.
\end{abstract}

\maketitle

\pagestyle{myheadings} \markboth{}{On a parabolic problem: Stability properties of Ground States}

\section{Introduction} In this paper we discuss the stability
properties of positive radial solutions for the following equation
\begin{equation}\label{laplace}
  \Delta u+ f(u,|x|)=0 \, ,
\end{equation}
which are positive steady states of the following Cauchy problem
\begin{align}
  &u_t= \Delta u+ f(u,|x|), \label{parab} \\ &u(x,0) =
  \phi(x), \label{data}
\end{align}
where $x \in \RR^n$, $n >2$, and $f=f(u, |x|)$ is a potential which is
null for $u=0$, superlinear in $u$, and supercritical in a sense that
will be specified just below.

Let $u(x,t;\phi)$ be the solution of \eqref{parab}--\eqref{data}.  The
analysis of the long time behavior of $u(x,t;\phi)$ is strongly based
on the separation properties of the radial solutions of
\eqref{laplace}.  If $u(x)$ is a radial solutions of (\ref{laplace}),
setting $U(r)=u(x)$ when $r=|x|$, we find that $U=U(r)$ solves
\begin{equation}\label{radsta}
  U''+\frac{n-1}{r} U'+ f(U,r)=0 \, ,
\end{equation}
where $``\,'\, "$ denotes the derivative with respect to $r$. In the
whole paper we denote by $U(r,\alpha)$ the unique solution of
\eqref{radsta} with the initial condition $U(0,\alpha)=\alpha>0$.

In the last decades the Cauchy problem \eqref{parab}--\eqref{data} has
raised a great interest, starting from the model case
$f(u,|x|)=u^{q-1}$, and it has been analyzed by several authors (see,
e.g., \cite{BF,FWY,GNW1,GNW2, BH, PY1,PY2,PY3,QS, W}).  Since in the
whole paper we are interested in positive solutions, there is no
ambiguity in using the notation $u^{q-1}$.  It is well known that the
behavior of solutions of \eqref{radsta}, and consequently of
\eqref{parab}, changes drastically as $q$ passes through some critical
values.  In this paper we focus on the case where $q>2^*:=
\frac{2n}{n-2}$, so that for any $\alpha>0$ the solution $U(r,\alpha)$
of \eqref{radsta} is positive and bounded for any $r>0$ (i.e. it is a
Ground State), and especially on the case $q \ge \sigma^*$, where
\begin{equation}\label{Josep}
  \sigma^*:= \left\{ \begin{array}{ll} \frac{(n-2)^2-4n+8
        \sqrt{n-1}}{(n-2)(n-10)} & \textrm{if $n>10$}, \\ +\infty &
      \textrm{if $n\le 10$},
    \end{array} \right.
\end{equation}
so that Ground States (GSs) gain some stability properties (see
\cite{W}). We recall that $2^*$ is the Sobolev critical exponent,
which is related to the compactness of the embedding of $L^q$ in
$H^1$, while $\sigma^*$ is the Joseph-Lundgren exponent, \cite{JL}.

When $2^*<q<\sigma^*$ all the GSs intersect each other indefinitely,
and this fact is used to construct suitable sub- and supersolution for
\eqref{laplace}. Then, it is possible to show that, in this range of
parameters, GSs determine the threshold between solutions of
\eqref{parab} that blow up in finite time, and solutions that exist
for any $t$ and fade away. 
\begin{theorem}\label{thmB}\cite{W,GNW1}
  Assume $f(u,r)=u^{q-1}$, $2^*<q< \sigma^*$.  Then
  \begin{itemize}
  \item[(1)] If there is $\alpha>0$ such that $\phi(x) \gneqq
    U(|x|,\alpha)$, then there is $T(\phi)$ such that $\lim_{t \to
      T(\phi)^-}\|u(t,x; \phi)\|_{\infty}=+\infty$.\\[-0.2 cm]
  \item[(2)] If there is $\alpha>\! 0$ such that $\phi(x) \lneqq
    U(|x|,\alpha)$, then $\lim_{t \to +\infty}\|u(t,x;
    \phi)\|_{\infty}=\!0$.
  \end{itemize}
\end{theorem}
On the other hand, when $q \ge \sigma^*$, GSs are well ordered, and
gain some stability properties as we will see just below.

In fact, already in \cite{W}, the whole argument was generalized to
embrace the so called Henon-equation, i.e. $f(u,r)=r^{\delta}
u^{q-1}$, where $\delta>-2$.  In this case there is a shift in the
critical exponents, so we find convenient to introduce the following
parameters (see Section~\ref{steadyasympt} below, see also \cite{BF}
for more details) which will be widely used through the whole paper:
\begin{equation}
  \label{eq:l0} l_s:=2 \frac{q+\delta}{2+\delta}\quad
  \textrm{ and }\quad m(l_s):=\frac{2}{l_s-2}= \frac{2+\delta}{q-2}.
\end{equation}
In this context, the previous discussion is still valid, but we have
stability whenever $l_s \ge \sigma^*$, and we lose it for
$2^*<l_s<\sigma^*$ (see \cite{W}). Notice that $l_s$ reduces to $q$
for $\delta=0$. In both cases the GSs, $U$, decay as $U(r) \sim
U(r,+\infty)=P_1 r^{-m(l_s)}$ for $r \to +\infty$, and $U(r,+\infty)$
is the unique singular solution of \eqref{radsta}.

To clarify the notion of stability we use, we need to introduce the
definitions of suitable weighted norms (see, e.g., \cite{GNW1}). We
set
\begin{equation}
  \begin{aligned}
    \|\psi\|_\lambda & := \underset{x\in \mathbb{R}^n} \sup |(1 +
    |x|^\lambda)\psi (x)|, \\ \vt{\, \psi \,}_{\lambda} & :=
    \underset{x\in \mathbb{R}^n} \sup \left| \frac{(1 +
        |x|^{\lambda)}}{[\ln(2+|x|)]}\psi (x)\right|,
  \end{aligned}
\end{equation}
where $\psi$ is continuous and $\la \in \RR$, $k \in \mathbb{N}$.
\begin{defin} We say that a GS, $U(|x|)=U(|x|, \alpha)$, is stable
  with respect to some norm $\|\cdot \|_\lambda$ if for every
  $\epsilon>0$ there exists $\delta >0$ such that, for $\| \varphi -
  U\|_\lambda< \delta$, we have $\|u(\cdot , t, \varphi) -
  U(|\cdot|)\|_\lambda < \epsilon$ for all $t> 0$.

  Further, we say that $U(|x|)$ is weakly asymptotically stable with
  respect to $\|\cdot \|_\lambda$ when $U(|x|)$ is stable with respect
  to $\|\cdot \|_\lambda$, and there exists $\delta >0$ such that
  $\|u(\cdot , t, \varphi) - U(|\cdot|)\|_{\lambda'} \to 0$ as $t\to
  \infty$, if $\|\varphi - U\|_\lambda< \delta$ (respectively, there
  exists $\delta >0$ such that $\|u(\cdot , t, \varphi) -
  U(|\cdot|)\|_{\lambda'} \to 0$ as $t\to \infty$ for all $\lambda' <
  \lambda$, if $\|\varphi - U\|_\lambda< \delta$).
\end{defin}
Let us consider the quadratic equation in $\la$
\begin{equation}\label{defLa}
  \la^2+ \Big(n-2-\frac{4}{q-2}\Big) \la +
  2\Big(n-2-\frac{2}{q-2}\Big)=0.
\end{equation}
Equation \eqref{defLa} admits two real and negative solutions, say
$\la_2 \le \la_1 <0$ if and only if $q \ge \sigma^*$, which coincide
if and only if $q=\sigma^*$. Gui et al in \cite{GNW1} proved the
following theorem.

\begin{theorem}\label{thmD}\cite{W,GNW1}
  Assume $f(u,r)=u^{q-1}$, $q \ge \sigma^*$.  Let $\la_2 \le \la_1$ be
  the roots of equation \eqref{defLa}.
  \begin{enumerate}
  \item If $q>\sigma^*$ any GS $U(r,\alpha)$ is stable with respect to
    the norm $\|\cdot \|_{m(q)+ |\la_1|}$ and weakly asymptotically
    stable with respect to the norm $\|\cdot \|_{m(q)+ |\la_2|}$. \\[-0.7 em]
  \item If $q=\sigma^*$ any GS $U(r,\alpha)$ is stable with respect to
    $\vt{\, \cdot\,}_{m(q)+ |\la_1|}$ and weakly asymptotically stable
    with respect to the norm $\|\, \cdot\,\|_{m(q)+ |\la_1|}$.
  \end{enumerate}
\end{theorem}

Actually, there is a number of results meant to extend the previous
analysis to more general potentials $f$ (see, e.g.,
\cite{Bae1,DLL,Dnew,YZ,BF}).  In particular the instability result
given by Theorem~\ref{thmB}, and the stability result
Theorems~\ref{thmD}, have been extended also to the following equation
\begin{equation} \label{eq:kmodel-1} u_t = \Delta u + k(r)r^\delta
  u^{q-1}, \quad \textrm{where $\delta>-2$}, \; \textrm{and $r=|x|$}
\end{equation}
assuming $k(r)$ decreasing, uniformly positive and bounded, in the
cases $l_s > \sigma^*$ (see \cite{DLL}), and $l_s=\sigma^*$ (see
\cite{Dnew}). In particular, these hypotheses implies that the
singular radial solution $U(r,+\infty)$ of \eqref{laplace} behaves
like $r^{-m(l_s)}$ both as $r \to 0$ and as $r \to +\infty$.

In such a case $q$ is replaced by $l_s$ and also the values of
$\la_1,\la_2$ change accordingly, i.e. they solve
\begin{equation}\label{defLaDelta}
  \la^2+ \Big(n-2-2\frac{2 + \delta}{q-2}\Big) \la+
  \frac{2+\delta}{q-2} \Big(n-2-\frac{2+\delta}{q-2}\Big)=0.
\end{equation}

In \cite{BF} we proposed a unifying approach which allows to extend
Theorem \ref{thmB} to a more general class of nonlinearities $f$,
including \eqref{eq:kmodel-1}, but also more involved dependence on
$u$.

The purpose of this paper is to continue the analysis of \cite{BF},
extending the stability results found in Theorem \ref{thmD} to a
larger class of $f$.  This purpose is achieved with an approach
obtained through the combination of the main ideas in \cite{W,GNW1,
  DLL}, techniques borrowed from the theory of non-autonomous
dynamical systems (see \cite{JPY,BF}), along with the use of some new
arguments.

As far as \eqref{eq:kmodel-1} is concerned we are able to drop the
assumption of boundedness on $k$ replacing it by the following:
\begin{equation} \label{asintotico} k(r) \sim r^{-\eta}, \, \textrm{
    as }\,\, r\to 0,\,\, \textrm{ with }\,\, 0\leq \eta < 2 + \delta.
\end{equation}
Then, we can allow two different behaviors for singular and slow decay
solutions (see \cite{BF}), namely: $U(r) \sim r^{-m(l_s)}$ as $r \to
+\infty$ and $U(r) \sim r^{-m(l_u)}$ as $r \to 0$, where
%
%
\begin{equation}\label{esempioM}
  l_u= 2\frac{q+\delta-\eta}{2+\delta-\eta}\quad \textrm{ and }\quad
  m(l_u)= \frac{2+\delta-\eta}{q-2}.
\end{equation}

So we prove the following.
\begin{thm} \label{main-1} Let $f(u,r)$ be as in \eqref{eq:kmodel-1},
  where $k(r)\in C^1$ satisfies \eqref{asintotico}, is decreasing, and
  $\underset{r\to +\infty}{\lim} k(r)>0$.  Then
  \begin{enumerate}
  \item If $l_s>\sigma^*$ any GS $U(r,\alpha)$ is stable with respect
    to the norm $\|\cdot \|_{m(l_s)+ |\la_1|}$ and weakly
    asymptotically stable with respect to the norm $\|\cdot
    \|_{m(l_s)+ |\la_2|}$.
  \item If $l_s=\sigma^*$ any GS $U(r,\alpha)$ is stable with respect
    to $\vt{\, \cdot\,}_{m(l_s)+ |\la_1|}$ and weakly asymptotically
    stable with respect to the norm $\|\, \cdot\,\|_{m(l_s)+ |\la_1|}$
  \end{enumerate}
\end{thm}

In fact, our approach is flexible enough to consider also a finite sum
of power in $u$, i.e.
\begin{align}
  & f(u, |x|) = k_1(|x|)r^{\delta_1}|u|^{q_1-1} +
  k_2(|x|)r^{\delta_2}|u|^{q_2-1}, 
  \label{eq:potential-1}
\end{align}
where $q_1<q_2$, $k_i=k_i(|x|)$, $i= 1,2$, are supposed to be $C^1$
(see Theorem~\ref{stabile} and Theorem~\ref{asint.stabile}, below).

Equation~\eqref{eq:potential-1} has been already considered by Yang
and Zhang in \cite{YZ}, but just in the particular situation of
$k_1(r)=k_2(r) \equiv 1$. We emphasize that, even if it is not stated
so clearly, in \cite{YZ} it is required that
\begin{equation}\label{condizione}
  \frac{2+\delta_2}{q_2-2}(q_2-q_1)+ \delta_1 < 0,
\end{equation}
which excludes the important case $\delta_1=\delta_2=0$. With these
assumptions, Yang and Zhang were able to prove
Theorem~\ref{thmD}--(1), replacing $q$ by $l_s= 2
\frac{q_2+\delta_2}{2+\delta_2}$, and changing the values of $m(l_s)$
and of $\la_i$ accordingly.

As a consequence of our main results we are able to generalize the
results in \cite{YZ} and to prove Theorem~\ref{thmD}, allowing $k_i$
to depend on $r$, and even to be unbounded, i.e.
\begin{equation} \label{asintotico1} k_1(r)\sim r^{-\eta_1}\,\,
  \textrm{ and }\,\, k_2(r)\sim r^{-\eta_2}, \,\, \textrm{ as }\,\,
  r\to 0,
\end{equation}
with $0\leq \eta_i < 2+ \delta_i$, $i=1, 2$. However we still need to
require \eqref{condizione}.

\begin{thm} \label{main-2} Let $f(u,r)$ be as in
  \eqref{eq:potential-1}, and assume \eqref{condizione}, and
  \eqref{asintotico1}. Suppose that both $k_1(r)r^{
    \frac{2+\delta_2}{q_2-2}(q_2-q_1)+\delta_1}$ and $k_2(r)$ are
  decreasing, $k_1(r)$ is positive and $k_2(r)$ is uniformly positive.
  Then, setting $l_s = 2 \frac{q_2 + \delta_2}{2 + \delta_2}$, we get
  the same conclusions as in Theorem~\ref{main-1}.
\end{thm}
Notice that we can deal with non-monotone functions $k_1(r)$.  Under,
these assumptions we are able to prove Theorem \ref{thmD}--(2) which
is new even in the case $k_1(r)=k_2(r) \equiv 1$ considered in
\cite{YZ}.

The main ingredients to obtain our results on \eqref{parab} are the
separation and the asymptotic properties of GSs.  The separation
properties are a result of independent interest, and generalize the
ones obtained in \cite[Theorems 1,2]{DLL2}, \cite[Theorem 2]{YZ2}. As
a consequence we also get Proposition \ref{ord3-}, which gives an
insight on the behavior of the singular solution of \eqref{radsta},
which seems to play a key role in determining the threshold between
blowing up and fading solutions (see the introduction in \cite{W}).

To prove weak asymptotic stability, we need a suitable asymptotic
expansion for GSs, which refines and generalizes the ones of
\cite{DLL,YZ} (see Proposition~\ref{sintetizzo}, below). In fact in
\cite{DLL,YZ} the highly nontrivial proof relies on an iterative
scheme developed by \cite{W} in a simpler  (but still nontrivial)
context. Here, we followed a different idea: in fact we have
proved an asymptotic results for nonlinear systems of ODEs, which
seems to be new to the best of our knowledge, and that, in our
opinion, is of intrinsic mathematical interest (even for systems of
ODEs).  In this more general framework the statement assumes a
more comprehensible aspect, and the proof is simplified, even if it is
still quite cumbersome;  We rely on the the appendix 
for a detailed proof of this lemma.

Now, we review briefly some results which have been proved just in the
setting of Theorems \ref{thmB}, \ref{thmD}.  First, using some sub-
and super-solutions constructed on the self-similar solutions,
\cite{GNW2,Na} proved that $U(|x|,\alpha)$ is weakly asymptotically
stable in the norm $\| \cdot \|_l$ for any
$m(q)+\la_1<l<m(q)+\la_2+2$.  Further Naito in \cite{Na} showed that
this result is optimal, i.e. in this range asymptotic stability does
not hold.  Moreover Gui et al. in \cite{GNW2} proved that GSs are not
even stable if we use too coarse, but surprisingly also too fine
norms, namely for $l<m(q)+\la_1$ and for $l \ge n$. Notice that we
have stability for $l=m(q)+\la_1$, but still there is a small gap for
$m(q)+\la_2+2 < n$. Similarly the null solution is weakly
asymptotically stable if $m(q) \le l < n$ and unstable otherwise,
\cite{GNW2}.

Moreover in a series of papers \cite{FWY2,HoYa, Na} the authors showed
that the speed at which the solution $u(t,x; \phi)$ converges depends
linearly on the weight used to measure the distance with respect to
the GS. Namely if $\|\phi(x)-U(|x|,\alpha)\|_l$ is small enough then
$t^{\nu} \|u(t,x;\phi)-U(|x|,\alpha)\|_{l'}$ is bounded for any $t>0$,
where $\nu=\frac{1}{2} \max\{l-l', l-m(q)-\la_1 \}$, whenever
$m(q)+\la_1<l<m(q)+\la_2+2$ and $0<l'<l$. The extension of these
results to more general non-linearities will possibly be the object of
future investigations.

To complete the picture we recall that, if either the assumptions of
Theorem~\ref{thmB} or of Theorem~\ref{thmD} are satisfied, following
\cite{BF} we can construct a family of subsolutions $\phi$ for
\eqref{radsta} with arbitrarily small $L^{\infty}$ norm and decaying
like $r^{2-n}$ for $r$ large, and such that the solution $u(t,x,\phi)$
blows up in finite time. This type of behavior contradicts the idea
that the decay of the singular solution, i.e $r^{-m(q)}$, is the
critical one to determine the threshold between fading and blowing up
solutions: The situation is more intricate. This results is in fact
extended to more general nonlinearities $f$, see \cite{BF}.



To conclude, we briefly recall that, when the non-linearity $f(u,r)$
becomes unbounded as $r \to 0$, in general it is not possible to find
classical solutions of \eqref{parab}--\eqref{data}.  However it is
still possible to obtain mild solutions assuming that $f(u,r)r^l$ is
bounded for $l>-2$, and in fact the solutions $u$ are classical for $x
\ne 0$ and $t>0$, and they are $C^{\alpha,\alpha/2}$ also for $x=0$
and $t=0$ for any $\alpha \in (0, l+2)$.  For an exhaustive exposition
about such a topic we refer to \cite{W} (see also \cite{BF}).

\smallskip \smallskip

\textbf{Plan of the paper.}  The paper is divided as follows: In
Section~\ref{steadyasympt} we collect all the preliminary results
concerning the solutions of \eqref{radsta}. We prove ordering
properties and asymptotic estimates for positive solutions of such a
problem.  Section~\ref{sec:stability} is devoted to the proof of the
main results of the paper (from which Theorems~\ref{main-1} and
\ref{main-2} follow directly).

\section{Ordering results and asymptotic estimates for the
stationary problem.}\label{steadyasympt} The results of this section,
which are crucial for our analysis, are obtained by applying the
Fowler transformation to \eqref{radsta}. For this purpose we need to
introduce some quantities that will appear frequently in the whole
paper, i.e.
\begin{eqnarray}\label{constants}
  m(l)=\frac{2}{l-2} \,, & A(l)=n-2-2m(l) \, , & B(l)=m(l)[n-2-m(l)],
\end{eqnarray}
where $l>2$ is a parameter (which is related to $l_s$ and $l_u$, in
\eqref{eq:l0} and in \eqref{esempioM}, respectively) whose role will
be explained few lines below.  Set
\begin{equation}\label{transf1}
  \begin{array}{c}
    r=\eu^s \, , \quad y_1(s,l)=U(\eu^s)\eu^{m(l) s} \,,\quad
    y_2(s,l)=\dot{y}_1(s,l) \\[0.3 em] g(y_1,s;l)=f(y_1 \eu^{-m(l) s},
    \eu^s) \eu^{(m(l)+2)s}
  \end{array}
\end{equation}
Throughout the paper ``$\,\dot{\,}\, $" will denote the
differentiation with respect to $s$ (recall that ``$\,'\,$" indicates
differentiation with respect to $r$).  Using these transformations we
pass from \eqref{radsta} to the following system:
\begin{equation} \label{si.na} \left( \begin{array}{c} \dot{y}_1 \\
      \dot{y}_2 \end{array}\right) = \left( \begin{array}{cc} 0 & 1 \\
      B(l) & -A(l)
    \end{array} \right)
  \left( \begin{array}{c} y_1 \\ y_2 \end{array}\right) -\left(
    \begin{array}{c} 0 \\ g(y_1,s;l)\end{array}\right).
\end{equation}
Here and in the sequel, we write $\bs{y}(s,\tau;
\bs{Q};\bar{l})=\big(y_1(s,\tau; \bs{Q};\bar{l}), y_2(s,\tau;
\bs{Q};\bar{l}) \big)$ to denote a trajectory of \eqref{si.na}, where
$l=\bar{l}$, evaluated at $s$ and departing from $\bs{Q} \in \RR^2$ at
$s=\tau$.

For illustrative purpose assume first $f(u,r)= r^{\delta} u^{q-1}$, so
we can set $\displaystyle l= 2\frac{q+\delta}{2+\delta}$ and
\eqref{si.na} reduces to the following autonomous system
\begin{equation} \label{si.a} \left( \begin{array}{c} \dot{y}_1 \\
      \dot{y}_2 \end{array}\right) = \left( \begin{array}{cc} 0 & 1 \\
      B(l) & -A(l)
    \end{array} \right)
  \left( \begin{array}{c} y_1 \\ y_2 \end{array}\right) -\left(
    \begin{array}{c} 0 \\ (y_1)^{q-1} \end{array}\right).
\end{equation}
In this case we passed from a singular non-autonomous ODE to an
autonomous system from which the singularity has been removed.  Also
note that when $\delta=0$ we can simply take $l=q$.  The sign of the
constants $A(l)$, $B(l)$ defined in \eqref{constants} determine
respectively if the system is sub- or supercritical, if there are slow
decay solutions ($B(l) \ge 0$) or if they do not exist ($B(l) < 0$).

\begin{remark} \label{slow-fast} We recall that, with the assumptions
  used in this paper, positive solutions $U(r)$ of \eqref{radsta} have
  two possible behaviour as $r \to 0$: \newline \emph{Regular},
  i.e. $\liro U(r) = \alpha>0$, or \emph{Singular}, i.e. $\liro U(r) =
  +\infty$. \\ Similarly as $r \to +\infty$ either $\lir U(r)r^{n-2}=
  \beta>0$ and we say that $U(r)$ has \emph{fast decay}, or $\lir
  U(r)r^{n-2}= +\infty$ and we say that $U(r)$ has \emph{slow decay}.

  In fact the behavior of singular and slow decay solutions can be
  specified better, see Proposition \ref{asympt} below), and
  Proposition \ref{sintetizzo}.
\end{remark}

\begin{figure}
  \centering \includegraphics*[totalheight=3.8cm]{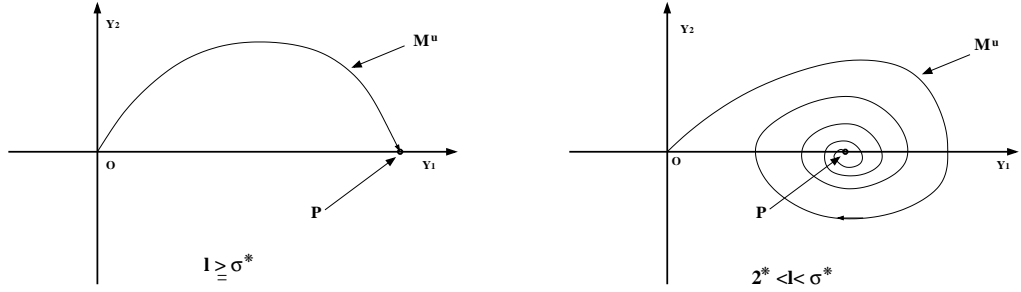}
  \caption{Sketches of the phase portrait of (\ref{si.na}), for $q>2$
    fixed.}\label{livelli}
\end{figure}

In this article we restrict the whole discussion to the case $l>2^*$,
therefore $A(l)>0$ and $B(l)>0$.  System (\ref{si.a}) admits three
critical points for $l>2^*$: The origin $O=(0,0)$, $\bs{P}=(P_1,0)$
and $-\bs{P}$, where $P_1= [B(l)]^{1/(q-2)}>0$.  The origin is a
saddle point and admits a one-dimensional $C^1$ stable manifold
$\overline{M}^s$ and a one-dimensional $C^1$ unstable manifold
$\overline{M}^u$, see Figure \ref{livelli}.  The origin splits
$\overline{M}^u$ in two relatively open components: We denote by $M^u$
the component which leaves the origin and enters the semi-plane $y_1
\ge 0$. Since we are just interested in positive solutions, with a
slight abuse of notation, we will refer to $M^u$ as the unstable
manifold.

\begin{remark}\label{criticalP}
  The critical point $\bs{P}$ of (\ref{si.a}) is a stable focus if
  $2^*<l<\sigma^*$ and a stable node if $l \ge \sigma^*$.
\end{remark}

As a consequence of some asymptotic estimates we deduce the following
useful fact (see, e.g. \cite{Fjdde, Fcamq}).

\begin{remark}\label{corrispondenze2}
  Let $u(r)$ be a solution of (\ref{radsta}) and let $\bs{Y}(s;l)$ be
  the corresponding trajectory of system (\ref{si.a}), with
  $l>2^*$. Then $u(r)$ is regular (respectively has fast decay) if and
  only if $\bs{Y}(s;l)$ converges to the origin as $s \to -\infty$
  (resp. as $s \to +\infty$), $u(r)$ is singular (respectively has
  slow decay) if and only if $\bs{Y}(s;l)$ converges to $\bs{P}$ as $s
  \to -\infty$ (resp. as $s \to +\infty$).
\end{remark}
Using the Pohozaev identity introduced in \cite{Po}, and adapted to
this context in \cite{FArch}, we can draw a picture of the phase
portrait of \eqref{si.a} (see Figure \ref{livelli} below) and deduce
information on positive solutions of (\ref{radsta}).  Then it is not
hard to classify positive solutions: In the supercritical case
($l>2^*$) all the regular solutions are GSs with slow decay, there is
a unique SGS with slow decay.

We stress that all the previous arguments concerning the autonomous
Equation \eqref{si.na} still hold true for any autonomous super-linear
system \eqref{si.na}. More precisely, whenever $g(y_1,s;l)\equiv
g(y_1;l)$ and $g(y_1;l)$ has the following property, denoted by
$\bs{G0}$ (see \cite{FAnn} for a proof in the general $p$-Laplace
context, see also \cite{BF}).
\begin{itemize}
\item[$\bs{G0}$:] There is $l>2^*$ such that $g(0;l)=0=\partial_{y_1}
  g(0,l)$ and $\partial_{y_1} g(y_1,l)$ is a positive strictly
  increasing function for $y_1>0$ and $\lim_{y_1 \to +\infty}
  \partial_{y_1} g(y_1,l)=+\infty$.
\end{itemize}
When $\bs{G0}$ holds true, we denote by $P_1$ the unique positive
solution in $y_1$ of $ g(y_1; l)=B(l) y_1$. Hence $(P_1,0)$ is again a
critical point for \eqref{si.a}.
Further, we let $\sigma_*<\sigma^*$ be the real solutions of the
equation in $l$ given by
\begin{equation}\label{defsigma}
  A(l)^2-4\big[\partial_y g(P_1,l)-B(l)\big]=0 ,
\end{equation}
which reduces to $A(l)^2-4(q-2)B(l)=0$ for $g(y_1)=(y_1)^{q-1}$.  We
emphasize that when $f(u,r)=u^{q-1}$ the value of $\sigma^*$ coincides
with the one given in (\ref{Josep}).  Notice that Remarks
\ref{criticalP}, \ref{corrispondenze2} continue to hold true in this
slightly more general context (see \cite{Fcamq,Fjdde}).

\subsection{ Main assumptions and
  preliminaries} \label{ssec:main-assumptions} We collect here below
the assumptions used in our main results:
\begin{itemize}
\item[$\bs{G1}$:] There exists $l_u \ge \sigma^*$ such that for any
  $y_1>0$ the function $ g(y_1,s;l_u)$ converges to a $s$-independent
  $C^1$ function $g(y_1,-\infty; l_u)\not\equiv 0$ as $s \to -\infty$,
  uniformly on compact intervals.  The function $g(y_1,s; l_u)$
  satisfies $\bs{G0}$ for any $s \in \RR$.  Further, there is
  $\varpi>0$ such that $\lim_{s \to -\infty} \eu^{-\varpi s}\partial_s
  g(y_1,s; l_u)=0$.  \\[-0.3 cm]
\item[$\bs{G2}$:] There exists $l_s \ge \sigma^*$ such that for any
  $y_1>0$ the function $ g(y_1,s; l_s)$ converges to a $s$-independent
  $C^1$ function $g(y_1,+\infty, l_s)\not\equiv 0$ as $s \to +\infty$,
  uniformly on compact intervals. The function $g(y_1, s; l_s)$
  satisfies $\bs{G0}$ for any $s \in \RR$.  Further, there is
  $\varpi>0$ such that $\lim_{s \to
    +\infty} \eu^{+\varpi s}\partial_s g(y_1,s;l_s)=0$.  \\[-0.3 cm]
\item[$\bs{G3}$:] Condition $\bs{G2}$ holds and $g(y_1,s;l_s)$ and
  $\partial_{y_1} g(y_1,s;l_s)$ are
  decreasing in $s$ for any $y_1>0$.  \\[-0.3 cm]
\item[$\bs{G4}$:] Condition $\bs{G2}$ is verified with $\varpi=\gamma$
  satisfying
  \[
  g(P_1^+,s;l_s)=g(P_1^+, {+\infty}; l_s)+ c \eu^{-\gamma s}
  +o(\eu^{-\gamma s})
  \]
  for a certain $c \ne 0$. \\[-0.3 cm]
\item[$\bs{K}$:] Either $f$ is as in \eqref{eq:kmodel-1} or $f$ is as
  in \eqref{eq:potential-1} and satisfies \eqref{condizione}.
\end{itemize}

Hypotheses $\bs{G1}$, $\bs{G2}$ are used to ensure that the phase
portrait of \eqref{si.na} converges to an autonomous system of the
form \eqref{si.a} (with $l \ge \sigma^*$), respectively as $s \to \pm
\infty$.  $\bs{G3}$ is needed to prove ordering properties of positive
solutions and generalizes the condition required in \cite{DLL}.
$\bs{G4}$ is needed to derive asymptotic estimates on slow decay
solutions of \eqref{radsta}, and it gives back the standard
requirement when $f(u,r)=k(r) u^{q-1}$, i.e. $k(r)=k(\infty)+c
r^{-\gamma}+o(r^{-\gamma})$ (see \cite{DLL}).  Actually, condition
$\bs{G4}$ is assumed for definiteness and may be weakened, at the
price of some additional cumbersome technicalities.  Finally,
condition \textbf{K} is a technical requirement we are not able to
avoid, which in fact is implicitly assumed also in \cite{YZ}. It
implies that there is $c>0$ such that
\begin{equation}\label{Kimplica}
  B(l_s)=  \frac{g(P_1^+, {+\infty}; l_s)}{P_1^+}= c   |P_1^+|^{\overline q-2} =
  \frac{\partial_{y_1}g(P_1^+, {+\infty}; l_s)}{  \overline q-1}
\end{equation}
with $\overline q = q$ in the case of \eqref{eq:kmodel-1}, and
$\overline q= q_2$ for the potential \eqref{eq:potential-1}.
\begin{remark}\label{spiegoluls}
  Observe that $\bs{G1}$ and $\bs{G2}$ are satisfied, e.g., in the
  following cases:
  \begin{itemize}
  \item{} For equation \eqref{eq:kmodel-1} with $k$ satisfying
    \eqref{asintotico}: $l_s$ and $l_u$ are as in \eqref{eq:l0}
    and \eqref{esempioM}, respectively. \\[-0.35 cm]
  \item{} When $f$ is as in \eqref{eq:potential-1} and
    \eqref{asintotico1} holds true: $l_s$ is as in Theorem
    \ref{main-2}, i.e.  $l_s= \min
    \big\{2\frac{q_i+\delta_i}{2+\delta_i} \mid i=1,2 \big\}$, while
    $l_u= \max \big\{2\frac{q_i+\delta_i-\eta_i}{2+\delta_i-\eta_i}
    \mid i=1,2 \big\}$. We also emphasize that, if we consider
    \eqref{eq:potential-1}, then \eqref{condizione} amounts to ask for
    $2\frac{q_2+\delta_2}{2+\delta_2} \le
    2\frac{q_1+\delta_1}{2+\delta_1}$; so $\bs{K}$ is not satisfied if
    $\delta_i=\eta_i=0$, since $l_s=q_1 < q_2=l_u$.
  \end{itemize}
\end{remark}

\begin{lemma}\label{decresce}
  Assume $\bs{G2}$ and $\bs{G3}$, then we have the following
  \begin{description}
  \item[$\bs{A^-}$] The function $G(y_1,s;2^*):=\int_0^{y_1}g(a,s;
    2^\ast) da$ is decreasing in $s$ for any $y_1 >0$ strictly for
    some $s$.
  \end{description}
\end{lemma}
\begin{proof}
  Set $G(z,s,l_s)= \int_0^z g(a,s,l_s) da$, $H(z,s)=G(z,s,l_s)/z$.
  Then
  \begin{equation*}\label{gcala}
    G(z,s,l_s)=\int_0^z \frac{g(a,s,l_s)}{a}a da \le \frac{g(z,s,l_s)}{z} \int_{0}^{z}a da= \frac{z g(z,s,l_s)}{2}
  \end{equation*}
  Therefore $zg-G \ge zg-2G \ge 0$.  Since $\partial_z H= (zg-G)/z^2$,
  then $H(z,s)$ is increasing in $z$ and decreasing in $s$ for
  $\bs{G3}$. Hence
  \begin{equation*}
    G(y_1,s,2^*)=G(y_1 \eu^{-\delta s},s,l_s) \eu^{\delta s}=H(y_1
    \eu^{-\delta s},s)y_1 \, ,
  \end{equation*}
  so we conclude that $G(y_1,s,2^*)$ is decreasing in $s$.
\end{proof}

Observe that $\bs{A^-}$ means that the system is supercritical with
respect to $2^*$, and this ensures the existence of GSs for
\eqref{radsta} (see e.g. \cite[Proposition 2.12]{BF}). In the sequel,
in some cases, it will be convenient to use the slightly weaker
condition $\bs{A^-}$, along with $\bs{G2}$, in place of the
combination of $\bs{G2}$ and $\bs{G3}$.

\subsection{The stationary problem: the spatial dependent case}
\label{subsec-stationary-prob}
Now we turn to consider (\ref{si.na}) in the $s$-dependent case. The
first step is to extend invariant manifold theory to the
non-autonomous setting.

Assume $\bs{G1}$. We introduce the following $3$-dimensional
autonomous system, obtained from (\ref{si.na}) by adding the extra
variable $z=\textrm{e}^{\varpi t}$, i.e.
\begin{equation}\label{si.naa}
  \begin{pmatrix}
    \dot{y}_{1} \\ \dot{y}_{2} \\ \dot{z}
  \end{pmatrix}
  = \begin{pmatrix} 0 & 1 &0 \\ B(l_u) & -A(l_u) & 0 \\ 0 & 0 & \varpi
  \end{pmatrix}
  \begin{pmatrix}
    y_1 \\ y_2 \\ z
  \end{pmatrix}-
  \begin{pmatrix}
    0 \\ g(y_1,\frac{\ln(z)}{\varpi};l_u)\\ 0
  \end{pmatrix}.
\end{equation}
Similarly if $\bs{G2}$ is satisfied we set $l=l_s$ and $\zeta(t)=
\eu^{-\varpi t}$ and we consider
\begin{equation}\label{si.naas}
  \begin{pmatrix}
    \dot{y}_{1} \\ \dot{y}_{2} \\ \dot{\zeta}
  \end{pmatrix} =
  \begin{pmatrix}
    0 & 1 &0 \\ B(l_s) & -A(l_s) & 0 \\ 0 & 0 & \varpi
  \end{pmatrix}
  \begin{pmatrix}
    y_1 \\ y_2 \\ \zeta
  \end{pmatrix}
  - \begin{pmatrix} 0 \\ g(y_1,-\frac{\ln(\zeta)}{\varpi};l_s)\\ 0
  \end{pmatrix}.
\end{equation}
The technical assumptions at the end of $\bs{G1}$, $\bs{G2}$ are
needed in order to ensure that the systems are smooth respectively for
$z=0$ and $\zeta=0$.

We recall that if a trajectory of \eqref{si.na} does not cross the
coordinate axes indefinitely then it is continuable for any $s \in
\RR$ (see e.g. \cite[Lemma 3.9]{Fcamq}, \cite{fDmF}).  Consider
\eqref{si.naa} {(respectively \eqref{si.naas}) each trajectory
  corresponding to a definitively positive solution $u(r)$ of
  \eqref{radsta} is such that its $\alpha$-limit set is contained in
  the $z=0$ plane (respectively its $\omega$-limit set is contained in
  the $\zeta=0$ plane). Moreover such a plane is invariant and the
  dynamics reduced to $z=0$ (respectively, $\zeta=0$) coincides with
  the one of the autonomous system (\ref{si.na}) where $g(y_1,s;
  l_u)\equiv g(y_1,-\infty; l_u)$ (respectively, $g(y_1,s;l_s)\equiv
  g(y_1,+\infty; l_s)$).

  Observe that the origin of \eqref{si.naa} admits a $2$-dimensional
  unstable manifold $\mathbf{W^u}(l_u)$ which is transversal to $z=0$
  (and a $1$-dimensional stable manifold $M^s$ contained in $z=0$).

  Following \cite{Fdie} (see also \cite{JPY}), for any $\tau \in \RR$
  we have that
  \begin{equation*}
    W^u(\tau;l_u)=\bs{W^u}(l_u) \cap \{ z= \eu^{\varpi \tau} \}\,\,
    \textrm{ and }\,\, W^u(-\infty;l_u)=\bs{W^u}(l_u) \cap \{ z= 0 \}
  \end{equation*}
  are 1-dimensional immersed manifolds, i.e. the graph of $C^1$
  regular curves.  Moreover, they inherit the same smoothness as
  (\ref{si.naa}) and (\ref{si.naas}), that is: Let $K$ be a segment
  which intersects $W^u(\tau_0;l_u)$ transversally in a point
  $\Q(\tau_0)$ for $\tau_0 \in[-\infty,+\infty)$, then there is a
  neighborhood $I$ of $\tau_0$ such that $W^u(\tau;l_u)$ intersects
  $K$ in a point $\Q(\tau)$ for any $\tau \in I$, and $\Q(\tau)$ is as
  smooth as (\ref{si.naa}).

  Since we need to compare $W^u(\tau;l_u)$ and $W^s(\tau;l_s)$, we
  introduce the manifolds:
  \begin{equation}\label{cambioL}
    W^u(\tau;l_s):= \big\{ \R=\Q
    \textrm{exp}\big(\big(m(l_s)-m(l_u)\big)\tau \big) \in \RR^2
    \mid \Q \in W^u(\tau;l_u) \big\}.
  \end{equation}
  Note that $W^u(\tau;l_u)$ and $W^u(\tau;l_s)$ are homothetic, since
  they are obtained from each other simply multiplying by an
  exponential scalar.  However, if $l_u >l_s$, $W^u(\tau;l_s)$ becomes
  unbounded as $\tau \to -\infty$.  In order to deal with bounded
  sets, we also define the following manifold which will be useful in
  Section~\ref{sec:stability}, i.e.
  \begin{align}\label{manistella}
    \begin{aligned} W^u(\tau;l_*):= & \left\{ \begin{array}{ll}
          W^u(\tau;l_u) & \textrm{if } \tau \le 0 \\ W^u(\tau;l_s) &
          \textrm{if } \tau \ge 0
        \end{array}
      \right.  , \quad \xi(\tau) := \left\{ \begin{array}{ll} z(\tau)
          & \textrm{if } \tau \le 0 \\ 2-\zeta(\tau) & \textrm{if }
          \tau \ge 0
        \end{array}   \right. \\
      \intertext{and} & \bs{W^u}(l_*):= \big\{ (\Q, \xi(\tau)) \mid \Q \in
      W^u(\tau;l_*) \big\}.
    \end{aligned} \hspace{-2 cm}
  \end{align}

  The sets $W^u(\tau;l_u)$ may be constructed also using the argument
  of \cite[\S 13]{CoLe}, simply requiring that \eqref{si.na} is $C^1$
  in $\bs{y}$ uniformly with respect to $t$ for $t \le \tau$ in a
  fixed neighborhood of the origin.  In this case we see that the
  tangent to $W^u(\tau;l_u)$ is the unstable space of the system
  obtained from \eqref{si.na} linearizing in the origin. So we get the
  following.

\begin{remark}\label{tangente}
  Assume $\bs{G1}$. Then, in the origin $W^u(\tau;l_u)$ is tangent to
  the line $y_2=m(l_u)y_1$, for any $\tau \in \RR$.  Since
  $W^u(\tau;l_u)$, $W^u(\tau;l_s)$ and $W^u(\tau;l_*)$ are homothetic,
  they are all tangent to $y_2=m(l_u)y_1$ in the origin.
\end{remark}

As in the $s$-independent case, we see that the \emph{regular
  solutions} correspond to the trajectories in $W^u$ (see \cite{Fdie,
  Fcamq}).  More precisely, from Lemma~3.5 in \cite{Fcamq}, we get the
following.

\begin{lemma}\label{corrispondenze}
  Assume $\bs{G1}$, $\bs{G2}$.  Consider the trajectory
  $\bs{y}(s,\tau,\Q;l_u)$ of (\ref{si.na}) with $l=l_u$, the
  corresponding trajectory $\bs{y}(t,\tau,\R; l_s)$ of (\ref{si.na})
  with $l=l_s$ and let $u(r)$ be the corresponding solution of
  (\ref{radsta}). Then $\R=\Q \textrm{exp}[ (m(l_s)-m(l_u)) \tau]$.

  Further $u(r)$ is a regular solution if and only if $\Q \in
  W^u(\tau; l_u)$ or equivalently $\R \in W^u(\tau; l_s)$.
\end{lemma}
Now, we turn to consider singular and slow decay solutions of
\eqref{radsta}.  Let $P_1^-$, $P_1^+$ be the unique positive solutions
in $y_1$ respectively of $B(l_u)y_1=g(y_1,-\infty;l_u)$ and of
$B(l_s)y_1=g(y_1,+\infty;l_s)$, and set $\bs{P^{\pm}}=(P_1^{\pm},0)$.
Then, it follows that $(\bs{P^-},0)$ and $(\bs{P^+},0)$ are
respectively critical points of (\ref{si.naa}) and \eqref{si.naas}.

If $l_u \ge 2^*$, then $(\bs{P^{-}},0)$ admits a $1$-dimensional
exponentially unstable manifold, transversal to $z=0$ (the graph of a
trajectory which will be denoted by $\bs{y^*}(s,*;l_u)$) for system
\eqref{si.naa}, while if $l_s > 2^*$ then $(\bs{P^{+}},0)$ is stable
for \eqref{si.naas}, so it admits a $3$-dimensional stable manifold
(an open set).

From \cite[Proposition 2.12]{BF} we find the following
\begin{prop}\label{super}\cite{BF}
  Assume $\bs{G1}$, $\bs{G2}$, and $\bs{A^-}$.  Then, all the regular
  solutions $U(r,\alpha)$ of (\ref{radsta}) are GSs with slow decay,
  there is a unique singular solution, say $U(r,\infty)$, and it is a
  SGS with slow decay.
\end{prop}
\begin{prop}\label{asympt}\cite{BF}
  Assume $\bs{G1}$, $\bs{G2}$.  Then if $u(r)$ and $v(r)$ are
  respectively a singular and a slow decay solution of (\ref{radsta})
  we have $u(r)r^{m(l_u)} \to P_1^-$ as $r \to 0$ and $u(r)r^{m(l_s)}
  \to P_1^+$ as $r \to +\infty$.
\end{prop}

\subsection{Separation properties of stationary solutions.}
In this section we adapt the argument of \cite{DLL} and of \cite{YZ}
to obtain separation properties of \eqref{radsta}.  We begin by the
following Lemma which is rephrased from \cite[Theorem 4.1]{YZ2}, which
is a slight adaption of \cite[Lemma 2.11]{DLL}. We emphasize that
condition $\bs{K}$ is needed to prove estimate \eqref{graph2} below,
and it is in fact implicitly required in \cite[Theorem 4.1]{YZ2}, even
if it is not explicitly stated.
\begin{lemma}\label{ord0}
  Assume $\bs{G1}$, $\bs{G2}$, $\bs{G3}$, $\bs{K}$. Let
  $\bs{\bar{y}}(s)$ be the trajectory of \eqref{si.na} corresponding
  to the GS $U(r,\alpha)$ of \eqref{radsta}.  Then, for any $s \in
  \RR$ we have $\bar{y}_2(s)=\dot{\bar{y}}_1(s) \ge 0$,
  $0<\bar{y}_1(s) < P^+_1$ and
  \begin{equation}
    g(\bar{y}_1(s),s;l_s) < B(l_s) \bar{y}_1(s) \label{eq:g}
  \end{equation}
\end{lemma}
\begin{proof}
  Let us recall that all the regular solutions are GSs, this is a
  direct consequence of Proposition~\ref{super} and
  Lemma~\ref{corrispondenze}.  Let $\bs{\bar{y}}(s;l_u)=
  \bs{\bar{y}}(s) \eu^{(\alpha_{l_s}-\alpha_{l_u})s}$ be the
  corresponding trajectory of \eqref{si.na} where $l=l_u$, then, by
  standard fact in dynamical system theory, see \cite{CoLe}, we see
  that there are $c_i>0$ such that $\bar{y}_i(s;l_u)
  \eu^{-\alpha_{l_u} s} \to c_i$ as $s \to -\infty$ for $i=1,2$. Hence
  $ \bar{y}_i (s) \sim c_i \eu^{\alpha_{l_s}s} \to 0$ as $s \to
  -\infty$ for $i=1,2$: So \eqref{eq:g} is satisfied for $s \ll 0$.

  Let us set
  \begin{equation}\label{s0}
    s_0:= \sup \big\{ S \in \RR \mid g(\bar{y}_1(s),s;l_s) <
    B(l_s)\bar{y}_1(s) \, \, \textrm{ for any $s<S$ } \big\},
  \end{equation}
  so that \eqref{eq:g} holds for $s<s_0$. \smallskip

  It follows that $\dot{\bar{y}}_2(s)+A(l_s)\bar{y}_2(s)>0$ for
  $s<s_0$, hence $w(s)= \bar{y}_2(s) \eu^{A(l_s) s}$ is increasing for
  $s<s_0$. Since $w(s)\to 0$ as $s \to -\infty$ we find that
  $\bar{y}_2(s)>0$, for $s \le s_0$.

  Further, assume by contradiction that there is $\tilde{s}< s_0$ such
  that $\bar{y}_1(\tilde{s})=P^+_1$. Then, from $\bs{G3}$, for $s<
  \tilde{s}$ we have
  \begin{equation*}
    g(\bar{y}_1(\tilde{s}),+\infty;l_s) \le
    g(\bar{y}_1(\tilde{s}),\tilde{s};l_s) <
    B(l_s)\bar{y}_1(\tilde{s})=g(P^+_1,+\infty;l_s).
  \end{equation*}
  Since $g(\cdot,+\infty;l_s)$ is increasing we get
  $\bar{y}_1(\tilde{s})<P_1^+$, and we have a contradiction. Thus,
  $0<\bar{y}_1(s)<P_1^+$ for $s<s_0$.

  Now, we show that $s_0=+\infty$, so that \eqref{eq:g} holds for any
  $s \in \RR$ and the Lemma is proved. Assume by contradiction that
  $s_0<+\infty$. Consider the curve $\bs{\bar{y}}(s) =(\bar{y}_1(s),
  \bar{y}_2(s))$ defined for $s \le s_0$.  Since
  $\bar{y}_2(s)=\dot{\bar{y}}_1(s)>0$ for $s \le s_0$, it follows that
  $\bs{\bar{y}}(s)$ is a graph on the $y_1$-axis, and we can
  parametrize it by $\bar{y}_1$. Hence, we set
  $Q(\bar{y}_1):=\dot{\bar{y}}_1(\bar{y}_1)$ so that $\bs{\bar{y}}(s)$
  for $s \le s_0$ and $\bs{\Gamma:=\Gamma}(y_1)=(y_1,Q(y_1))$ for $y_1
  \in (0,\bar{y}_1(s_0)]$ represent the same curve.  As a consequence
  we have that
  \begin{equation}\label{graph}
    \frac{\partial Q}{\partial \bar{y}_1}=\frac{\partial Q}{\partial s}
    \frac{\partial s}{\partial
      \bar{y}_1}=\frac{\ddot{\bar{y}}_1}{\dot{\bar{y}}_1}=
    -A(l_s)+\frac{B(l_s)\bar{y}_1-g(\bar{y}_1,s;l_s)}{Q(\bar{y}_1)}.
  \end{equation}
  In the phase plane, consider the line $r(\mu)$ passing through
  $\bs{R}=(\bar{y}_1(s_0),0)$ with angular coefficient $-\mu$, i.e.
  \begin{equation*}
    r(\mu):= \big\{ (y_1,y_2) \mid y_2= \mu
    (\bar{y}_1(s_0)-y_1) \big\}.
  \end{equation*}
  Since $\bar{y}_2(s_0)=\dot{\bar{y}}_1(s_0)>0$, we see that
  $\bs{\Gamma}(\bar{y}_1(s_0))=(\bar{y}_1(s_0), \bar{y}_2(s_0))$ lies
  above $\bs{R}$.  By construction $r(\mu)$ intersects $\bs{\Gamma}$
  at least in a point, for any $\mu>0$: We denote by
  $\big(Y_1(\mu),\mu (\bar{y}_1(s_0)-Y_1(\mu))\big)$ the intersection
  with the smallest $Y_1$. Then, it follows that $Y_1<\bar{y}_1(s_0)$
  and $\frac{\partial Q}{\partial \bar{y}_1}(Y_1) \ge -\mu$. From
  these inequalities, along with \eqref{graph}, and using the fact
  that
  \begin{equation} \label{eq:rapporto-B} B(l_s)\bar{y}_1(s_0)=
    g(\bar{y}_1(s_0),s_0;l_s)
  \end{equation}
  we get
  \begin{equation}\label{graph2}
    \begin{aligned}
      -\mu \le \! \frac{\partial Q}{\partial \bar{y}_1}(Y_1) & =\!
      -A(l_s)+\frac{ B(l_s) [Y_1- \bar{y}_1(s_0)]\! + \! [
        g(\bar{y}_1(s_0),s_0;l_s) -g(Y_1,s;l_s)]}{\mu
        [\bar{y}_1(s_0)-Y_1]} \\ &\le \!
      -A(l_s)-\frac{B(l_s)}{\mu}+\frac{g(\bar{y}_1(s_0),s_0;l_s)
        -g(Y_1,s_0;l_s)}{\mu [\bar{y}_1(s_0)-Y_1]} \\ & \le \!
      -A(l_s)+\frac{1}{\mu } \left[ - B +\partial_{y_1} g(C,s_0;l_s) \right]  \\
      & \le \! -A+\frac{1}{\mu }  \left[ - B +\frac{(\bar{q}-1) g(C,s_0;l_s)}{C} \right] \\
      & < \! -A+\frac{1}{\mu } \left[ - B +\frac{(\bar{q}-1)
          g(\bar{y}_1(s_0),s_0;l_s)}{\bar{y}_1(s_0)} \right] =
      -A+\frac{B(\bar{q}-2)}{\mu }
    \end{aligned}\hspace{-0.38 cm}
  \end{equation}
  where $C \in (Y_1,\bar{y}_1(s_0))$ and we used the mean value
  theorem.  Further $\bar{q}$ stands for $q$ if $f$ is of type
  \eqref{eq:kmodel-1} and it stands for $q_2$ if $f$ is of type
  \eqref{eq:potential-1}.  Therefore, using \eqref{graph2} along with
  \eqref{eq:rapporto-B}, we obtain
  \begin{equation*}
    \mu^2 - A \mu +B(\bar{q}-2)= \mu^2 - A \mu- B+ \partial_{y_1}g(P_1^+,+\infty, l_s)
    > 0, \, \,\textrm{ for any $\mu>0$.}
  \end{equation*}
  But this is verified if and only if
  \begin{equation*}
    A^2 - 4   B(\bar{q}-2) =  A^2-4[ \partial_{y_1}g(P_1^+,+\infty, l_s)-B]
    < 0 ,
  \end{equation*}
  which is equivalent to $l_s \in (\sigma_*,\sigma^*)$, cf
  \eqref{defsigma}, so we have found a contradiction. Hence
  $s_0=+\infty$. In particular, it follows that
  $\bar{y}_1(s)<P_1^{+}$, $\dot{\bar{y}}_1(s)>0$, for any $s \in \RR$,
  and \eqref{eq:g} holds true.
\end{proof}

\begin{remark}\label{ord1}
  Assume the hypotheses of Lemma \ref{ord0}; Assume further $\bs{K}$,
  then
  \begin{equation} \label{eq:g1} \frac{\partial{g}}{\partial
      y_1}(\bar{y}_1(s),s;l_s) < \frac{\partial{g}}{\partial
      y_1}(P_1^+,+\infty;l_s)
  \end{equation}
\end{remark}
\begin{proof}
  From a straightforward computation we see that, when $f$ is as in
  \eqref{eq:kmodel-1}, then \eqref{eq:g} implies \eqref{eq:g1}. When
  $f$ is as in \eqref{eq:potential-1} then
  \begin{equation*}
    \partial_{y_1}g(y_1,s,l_s)=(q_1-1) k_1(\eu^s) y_1^{q_1-2} + (q_2-1)
    k_2(\eu^s) y_1^{q_2-2} \le (q_2-1) g(y_1,s,l_s)/y_1.
  \end{equation*}
  So, let $\bs{\bar{y}}(s)$ be a trajectory corresponding to a GS of
  \eqref{radsta} as above; If $\bs{K}$ holds, from \eqref{eq:g} we get
  \begin{equation*}
    \frac{\partial g}{\partial y_1}(\bar{y}_1,s,l_s) \le (q_2-1) \frac{
      g(\bar{y}_1(s),s,l_s)}{\bar{y}_1(s)} \le (q_2-1) \frac{
      g(P_1^+,+\infty;l_s)}{P_1^+} \le \frac{\partial g}{\partial
      y_1}(P_1^+,+\infty;l_s),
  \end{equation*}
  so \eqref{eq:g1} follows and the Lemma is proved.
\end{proof}

\begin{prop}\label{ord1-}
  Assume $\bs{G1}$, $\bs{G2}$, $\bs{G3}$. Assume further
  $\bs{K}$. Then $U(r,\alpha_1)<U(r,\alpha_2)$ for any $r>0$, whenever
  $\alpha_1<\alpha_2$.
\end{prop}

We emphasize that if $g(y_1,s;l)$ is $s$-independent, as in \cite{W},
Lemma \ref{ord0} implies Proposition~\ref{ord1-}. This fact follows
directly by noticing that $M^u$ is a graph on the $y_1$-axis, since
$y_1(s)=U(\eu^s,\alpha) \eu^{m(l_s)s}$ is increasing in $s$, for any
$\alpha>0$.  In view of Lemma~\ref{ord0}, we can parametrize the
manifold $M^u$ by $\alpha$, then the ordering of the regular solutions
$U(r,\alpha)$ is preserved as $s$ varies (i.e. as $r$ varies), since
they all move along a $1$-dimensional object.

When we turn to consider an $s$-dependent function $g(y_1,s; l)$,
Proposition \ref{ord1-} needs a separate proof, which can be obtained
by adapting the ideas developed in \cite{DLL,YZ}.  In fact, in such a
case $W^u(\tau;l_s)$ is still one dimensional but may not be a graph
on the $y_1$-axis, so a priori we may lose the ordering property.

\begin{proof}[Proof of Proposition~\ref{ord1-}]
  Let us set $Q(s)=\eu^{\la_1 s}$ and observe that
  \begin{equation}\label{linear}
    \ddot{Q}+A\dot{Q}+[\partial_{y_1}g (P_1^+,+\infty;l_s)-B] Q=0.
  \end{equation}
  Denote by $W(s):= [U(\eu^{s},\alpha_2)- U(\eu^{s},\alpha_1)]
  \eu^{m(l_s) s }$, and observe that
  \begin{equation}\label{difference}
    \ddot{W}+A\dot{W}-B W+ D(s)=0,
  \end{equation}
  where $D(s):=g(U(\eu^{s},\alpha_2)] \eu^{m(l_s) s }, s;l_s)-
  g(U(\eu^{s},\alpha_1)] \eu^{m(l_s) s }, s;l_s)$.

  Using continuous dependence on initial data we see that
  $U(r,\alpha_2)>U(r,\alpha_1)$ for $r$ small enough, so that $D(s)>0$
  for $s \ll 0$.  Assume by contradiction that there is
  $\bar{r}=\eu^{\bar{s}}>0$ such that $U(r,\alpha_2)- U(r,\alpha_1)>0$
  for $0\le r<\bar{r}$, and $U(\bar r, \alpha_2)- U(\bar r,
  \alpha_1)=0$.  Then, $W(s)$, and $D(s)$ are positive for $s<\bar{s}$
  and they are null for $s=\bar{s}$.

  Setting $Z(s) :=\dot{W}(s)Q(s)-W(s)\dot{Q}(s)$, by direct
  calculation we can easily see that
  $\dot{Z}(s)=\ddot{W}(s)Q(s)-W(s)\ddot{Q}(s)$. Then from
  \eqref{linear} and \eqref{difference} we get
  \begin{equation}\label{ODE}
    \dot{Z}=-A Z(s)+Q(s)[\partial_{y_1}g
    (P_1^+,+\infty;l_s)W(s)-D(s)].
  \end{equation}
  Observe now that $W(s) \sim (\alpha_2-\alpha_1) \eu^{m(l_s)s}$, as
  $s \to -\infty$, and also that
  \begin{equation*}
    \dot{W}(s)=m(l_s)W(s)+
    [U'(\eu^{s},\alpha_2)-U'(\eu^{s},\alpha_1)]\eu^{[1+m(l_s)]s}\sim
    m(l_s) (\alpha_2-\alpha_1) \eu^{m(l_s)s},
  \end{equation*}
  as $s \to -\infty$. Hence, we get
  \begin{equation}\label{zero}
    Z(s)\sim
    (m(l_s)-\la_1(l_s))(\alpha_2-\alpha_1)(\eu^{[m(l_s)+\la_1(l_s)]s})
    \to 0 \quad \textrm{ as $s \to -\infty$ .}
  \end{equation}
  Moreover $\la_1(l_s)+A(l_s)=-\la_2(l_s)>0$ and $D(s) \to 0$ as $s
  \to -\infty$, hence $\eu^{A s} Q(s)D(s)\in
  L^1(-\infty,\bar{s}]$. Since $Z(s)$ is the unique solution of
  \eqref{ODE} satisfying \eqref{zero} we find
  \begin{equation}\label{integral}
    Z(\bar{s})= \int_{-\infty}^{\bar{s}} \eu^{-A (\bar{s}-s)}
    Q(s)[\partial_{y_1}g (P_1^+,+\infty;l_s)W(s)-D(s)] ds.
  \end{equation}
  From the mean value theorem we find
  \begin{equation*}
    \partial_{y_1}g (P_1^+,+\infty;l_s)W(s)-D(s)= [\partial_{y_1}g
    (P_1^+,+\infty;l_s)- \partial_{y_1}g (U(s),s;l_s)]W(s).
  \end{equation*}
  where $U(s)$ lies between $U(r,\alpha_1)r^{m(l_s)}$ and
  $U(r,\alpha_2)r^{m(l_s)}$. Since $\partial_{y_1}g (y_1,s;l_s)$ is
  increasing in $y_1$, and using \eqref{eq:g1}, for $s <\bar{s}$ we
  find
  \begin{equation}\label{convex}
    \begin{split}
      \partial_{y_1}g (P_1^+, &\,+\infty;l_s)W(s)-D(s) \\ & \ge
      [\partial_{y_1}g (P_1^+,+\infty;l_s)- \partial_{y_1}g
      (U(\eu^s,\alpha_2)\eu^{m(l_s)s},s;l_s)]W(s)>0
    \end{split}
  \end{equation}

  Hence, from \eqref{integral} and \eqref{convex} we get
  \begin{equation*} 
    0<Z(\bar{s})=\dot{W}(\bar{s})Q(\bar{s})-W(\bar{s})\dot{Q}(\bar{s})=
    \dot{W}(\bar{s})Q(\bar{s}) \, ,
  \end{equation*}
  which gives $\dot{W}(\bar{s})>0$. Thus, we find $W(s)<0$ for
  $|s-\bar{s}|$ small enough, and this gives a contradiction.
  Therefore, $U(r,\alpha_2)- U(r,\alpha_1)>0$ for any $ r \ge 0$.
\end{proof}

Now, we consider the singular solution $U(r,\infty)$.
\begin{prop}\label{ord3-}
  Assume the hypotheses of Proposition \ref{ord1-}, then
  $U(r,\infty)r^{m(l_s)}$ is non-decreasing for any $r>0$, and
  $U(r,\alpha)<U(r,\infty)$ for any $r>0$, $\alpha>0$.
\end{prop}
Actually, this result is new even for $f$ of both types $f(u,r)=K(r)
u^{q-1}$ and $f(u,r)= u^{q_1-1}+u^{q_2-1}$, which are considered in
\cite{DLL} and \cite{YZ}, respectively.
\begin{proof}
  The result is well known when the system is autonomous: In fact in
  this case $U(r, \infty) r^{m(l_s)}\equiv P_1^{+}$ and
  $W^u_{l_s}=W^u_{l_u}$ is a graph on the $y_1$-axis connecting the
  origin and $\bs{P^{+}}$.

  Now, we turn to consider the $s$-dependent setting. From the
  previous discussion we know that the manifold $M^u$ of the
  autonomous system \eqref{si.na}, where $l=l_u$,
  $g=g(y_1,-\infty;l_u)$, is a graph on the $y_1$-axis connecting the
  origin and the critical point $\bs{P^{-}}$. Moreover, observe that
  for any $\tau \in \RR$ the manifold $W^u(\tau;l_u)$ is a graph
  connecting the origin and the unique trajectory $\bs{y^*}(s; l_u)$
  $U(r,\infty)$ (and such that $\lim_{s \to -\infty}\bs{y^*}(s;
  l_u)=\bs{P^{-}}$).

  We claim that $W^u(\tau;l_u)$ is a graph on the $y_1$-axis, for any
  $\tau\in \RR$.  In fact let $\Q, \R \in W^u(\tau;l_u)$, with
  $\Q=(Q_1,Q_2), \R=(R_1,R_2)$, and let $U(r, \alpha_Q)$ and $U(r,
  \alpha_R)$ be the corresponding solution of (\ref{radsta}).  From
  Proposition \ref{ord1-} we know that if $\alpha_Q<\alpha_R$, then
  \begin{equation}\label{abo}
    Q_1=U(\eu^{\tau}, \alpha_Q)\eu^{m(l_u)\tau}<U(\eu^{\tau},
    \alpha_R)\eu^{m(l_u)\tau}=R_1 \, ,
  \end{equation}
  so the claim follows.

  Moreover, we also get $Q_1< y_1^*(\tau;l_u)$. In fact assume by
  contradiction that $Q_1 > y_1^*(\tau;l_u)$. Then we can choose $\R$
  in the branch of $W^u(\tau;l_u)$ between $\Q$ and $\bs{y^*}(\tau;
  l_u)$, so that $\alpha_{R}>\alpha_{Q}$ and $Q_1 >R_1> y_1^*(\tau;
  l_u)$; but this contradicts \eqref{abo}.  Similarly if $Q_1 =
  y_1(\tau,*; l_u)$, then $\R \in W^u(\tau;l_u)$ is such that
  $\alpha_R>\alpha_Q$, and $R_1>Q_1= y_1(\tau,*; l_u)$.  But again we
  can choose $\bs{\tilde{R}}$ in the branch of $W^u(\tau;l_u)$ between
  $\R$ and $\bs{y^*}(\tau; l_u)$, and reasoning as above we find again
  a contradiction.  Therefore $U(r,\alpha)<U(r,\infty)$ for any $r>0$,
  and any $\alpha>0$.

  Further, since $W^u(\tau;l_u)$ and $W^u(\tau;l_s)$ are homothetic,
  cf \eqref{cambioL}, then $W^u(\tau;l_s)$ is a graph on the
  $y_1$-axis, which connects the origin and the trajectory
  $\bs{y^*}(s; l_s)$ corresponding to $U(r,\infty)$. Further
  $W^u(\tau;l_s) \subset \{(y_1,y_2) \mid 0<y_1<P_1^+ , \; y_2>0 \}$
  (see Lemma~\ref{ord0}).  Therefore $y_2^*(s; l_s) \ge 0$ for any $s
  \in \RR$.  Hence $U(r,\infty)r^{m(l_s)}$ is non-decreasing for any
  $r>0$, and the proof is concluded.
\end{proof}
Proposition~\ref{ord3-} is interpreted as follows in terms of system
\eqref{si.na}. 

\begin{remark}\label{ord3-bis}
  In the hypotheses of Proposition \ref{ord1-}, hence of
  Proposition~\ref{ord3-}, we have that $W^u(\tau;l_u)$,
  $W^u(\tau;l_s)$, and $W^u(\tau;l_*)$ are graphs on the $y_1$-axis
  respectively for any $\tau \in \RR$. Further they are contained in
  $y_2 \ge 0$ and connect the origin respectively with
  $\bs{y^*}(\tau;l_u)$, $\bs{y^*}(\tau;l_s)$, and
  $\bs{y^*}(\tau;l_*)$.
\end{remark}

 \subsection{Asymptotic estimates for slow decay solutions.}\label{subsec:asympotic}
 In this subsection we state the asymptotic estimates for slow decay
 solutions of \eqref{radsta}, which are crucial to prove our main
 results: We always assume $\bs{G1}$, $\bs{G2}$, and $\bs{G4}$.

 In fact, we generalize the results obtained for $f(u,r)=k(r)u^{q-1}$
 in \cite[ \S 3]{DLL} for $q>\sigma^*$, and in \cite{Dnew} for
 $q=\sigma^*$.  The main argument in \cite{DLL} has been re-used in
 \cite{YZ}, and it is an adaptation to the non-autonomous context of
 the scheme introduced by Li in \cite{L} (and developed in
 \cite{GNW1}).  Here, we follow a different approach, so we give an
 interpretation in terms of general facts in ODE theory of the
 argument behind the whole \cite[\S 3]{DLL}, which thereafter becomes
 clearer in our opinion.

 Due to assumption $\bs{G4}$ we can now set $\zeta= \eu^{-\gamma s}$
 in (\ref{si.naas}), and obtain a smooth system which has
 $\bs{\mathcal{P}}:=(P_1^+, 0,0)$ as critical point.  In this
 subsection we consider this system and its linearization around
 $\bs{\mathcal{P}}$ so we leave the explicit dependence on $l_s$
 unsaid.
 \begin{equation}\label{si.lin}
   \begin{pmatrix}
     \dot{y}_{1} \\ \dot{y}_{2} \\ \dot{\zeta}
   \end{pmatrix} =
   \begin{pmatrix}
     0 & 1 & 0 \\ B- \partial_{y_1} g^{+\infty}(P_1^+) & -A & 0 \\ 0 &
     0 & -\gamma
   \end{pmatrix}
   \begin{pmatrix}
     y_1 \\ y_2 \\ \zeta
   \end{pmatrix}
 \end{equation}
 Let us denote by $\mathcal{A}$ the matrix in \eqref{si.lin}: It has
 $3$ negative eigenvalues $\la_2 \le \la_1 <0$ and $-\gamma<0$
 ($\bs{G4}$ is needed in order to guarantee smoothness of the system
 \eqref{si.naas} for $\zeta=0$).  Therefore the critical point
 $\bs{\mathcal{P}}$ of \eqref{si.naas} is a stable node.

 Assume first that the $3$ eigenvalues are simple, then we have $3$
 eigenvectors, respectively $v_{1}=(1,-m+\la_1,0)$,
 $v_2=(1,-m+\la_2,0)$, and $v_z:=v_3=(0,0,1)$.  Any solution $\ell(t)$
 of \eqref{si.lin} can be written as
 \begin{equation}\label{lin}
   \ell(s)=\bar{a} v_1 \eu^{\la_1 s}+ \bar{b} v_2 \eu^{\la_2 s}+z v_z
   \eu^{-\gamma s}
 \end{equation}
 for some $\bar{a},\bar{b},z \in \RR$.

 By standard facts in invariant manifold theory (see \cite[\S
 13]{CoLe}), any trajectory $(\bs{y}(s), \zeta(s))$ of \eqref{si.naas}
 converging to $\bs{\mathcal{P}}$ can be seen as a non-linear
 perturbation of a solution $\ell(s)$ of \eqref{si.lin}. More
 precisely set $\bs{n}(s)=(n_1(s),n_2(s))= (y_1(s)- P_1^+,y_2(s))$,
 then $\bs{N}(s):=(n_1(s),n_2(s),\zeta(s))=\ell(s)+O(|\ell(s)|^2)$.
 Therefore
 \[
 n_1(s)= \bar{a} \eu^{\la_1 s}+ \bar{b} \eu^{\la_2 s}+z \eu^{-\gamma
   s} + O(\eu^{2\la_1 s}+ \eu^{2\la_2 s}+ \eu^{-2\gamma s})
 \]
 In the  appendix we prove that the expansion can be
 continued to an arbitrarily large order: This is the contained of
 Proposition \ref{sintetizzo2} and of its general form containing
 resonances, i.e.  Proposition \ref{sintetizzo}.  Let us rewrite
 \eqref{si.na} as
 \begin{equation}\label{xxx}
   \dot{\vec{x}}= \mathcal{A} \vec{x}+ \vec{N}(\vec{x})
 \end{equation}
 where $\vec{x}=(y_1,y_2,\zeta)$, and $\mathcal{A}$ is as the matrix
 in \eqref{si.lin}.
 \begin{prop}\label{sintetizzo2}
   Assume for simplicity $\vec{N} \in C^{\infty}$ and that the
   eigenvalues of $\mathcal{A}$ are real, negative and simple and are
   rationally independent, i.e there is no
   $\chi=(\chi_1,\chi_2,\chi_3) \in \mathbb{Z}^3 \backslash \{(0,0,0)
   \}$ such that $ \chi_1 |\la_1|+\chi_2 |\la_2|+\chi_3 \gamma=0$, so
   that no resonances are possible.  Further assume for definiteness
   that $|\la_1|<\gamma$.

   Then for any $k \in \mathbb{N}$ we can find a polynomial $P$ of
   degree $k$ in $3$ variables such that
  $$ y_1(t)= P(\eu^{\la_1 t}, \eu^{\la_2 t}, \eu^{-\gamma t})+o(\eu^{[(k+1)\la_1+\ep] t}) $$
  as $t \to +\infty$, for $\ep>0$ small enough.
\end{prop}
We remand the interested reader to  the Appendix for
details.

Now we state the result in a form which is more suitable for our
purpose; Set
\begin{equation}\label{Itheta}
  I_\theta=\big\{\chi=(\chi_1,\chi_2,\chi_3) \in \mathbb{N}^3 \, \colon \,  \chi_1
  |\la_1|+\chi_2 |\la_2|+\chi_3 |\gamma| \le \theta \big\}.
\end{equation}
Then, we can expand $n_1(s)$ as follows
\begin{equation}\label{expand0}
  \begin{array}{c}
    \ds n_1(s)= a \eu^{\la_1 s}+ b \eu^{\la_2 s}+z \eu^{-\gamma s} +
    P_{\theta}(s) +o( \eu^{-\theta s}),
  \end{array}
\end{equation}
where the function $P_{\theta}(s)$ is completely determined by the
values of the coefficients $a,b,z$.

As a first case, assume that $\gamma, |\la_1|, |\la_2|$ are rationally
independent. Then, there are constants $c^\chi \in \RR$ such that
\begin{equation}\label{Ptheta}
  \ds P_{\theta}(s)=\sum_{ \tiny \begin{array}{c}
      \chi\in I_\theta ,\\   \!\!|\chi| \ge 2 \end{array} }
  \!\!\!\! c^\chi \eu^{(\chi_1 \la_1+\chi_2 \la_2-\chi_3 \gamma)s} \,\,
  \textrm{ with $\chi=(\chi_1,\chi_2,\chi_3)$}
\end{equation}
and $|\chi|=\chi_1+\chi_2+\chi_3$.

Let us now consider the resonant cases, i.e.  when there are $M^0,
M^1, \ldots , M^j$, (a $j$-ple resonance)
$M^i=(\chi^i_1,\chi^i_2,\chi^i_3) \in I_\theta$, $|M^i|>0$ for $i=1,
\ldots, j$, such that
\begin{equation*}
  \chi^i_1 |\la_1| + \chi^i_2 |\la_2| + \chi^i_3 \gamma= \bar{\theta} \le
  \theta \, \;
\end{equation*}
Then, we have to replace $\sum_{i=0}^j c_{M^i} \eu^{(\chi^i_1
  \la_1+\chi^i_2 \la_2-\chi^i_3 \gamma)s}$ by
\begin{equation}\label{sost}
  \sum_{i=0}^j c_{M^i} s^{i} \eu^{(\chi^i_1 \la_1+\chi^i_2 \la_2-\chi^i_3
    \gamma)s} \; \textrm{ in the function $P_{\theta}$},
\end{equation}
(notice that we have included the possible case of resonances with the
linear terms, e.g., $\chi_2$ multiple of $\chi_1$ etc...).  The same
happens when we have resonances within the linear terms,
e.g. $|\la_1|=|\la_2|$ (i.e. $l_s=\sigma^*$), or $|\la_1|=\gamma$: We
replace the terms as done in \eqref{sost}.

Before collecting all these facts in Proposition \ref{sintetizzo}
below, we need some further notation.  Let us introduce the following
sets, i.e.
\begin{eqnarray}
  J_{|\la_1|} &=& \{\chi=(0,0,\chi_3) \in \mathbb{N}^3 \mid 0< \chi_3 \gamma < |\la_1| \}
  \, , \\ J_{|\la_2|} &=& \{\chi=(\chi_1,\chi_2,\chi_3) \in \mathbb{N}^3
  \mid |\la_1| < \chi_1|\la_1|+ \chi_3 \gamma < |\la_2| \}.
\end{eqnarray}
Observe that $ J_{|\la_1|}$ is empty if $|\la_1| \le \gamma$, and $
J_{|\la_2|}$ is empty, e.g., if $|\la_2| < 2|\la_1|$, and $|\la_2|\le
\gamma$.  We denote by
\begin{eqnarray}
  \ds \Psi(s) &=&\sum_{\chi=(0,0,\chi_3) \in J_{|\la_1|}} c^\chi \eu^{-\chi_3
    \gamma s} + \chi_r(s) \eu^{\la_1 s} \label{defPsi}
\end{eqnarray}
where $\chi_r(s)=0$ if $|\la_1|/\gamma \not\in \mathbb{N}$, and
$\chi_r(s)=\chi_r s$ if $|\la_1|/\gamma \in \mathbb{N}$ and
$l_s>\sigma^*$, while $ \chi_r(s)=\chi_r s^2$ if $|\la_1|/\gamma \in
\mathbb{N}$ and $l_s=\sigma^*$, for a certain $\chi_r \in \RR$.
\smallskip

  \begin{prop}\label{sintetizzo}
    Assume $\bs{G1}$, $\bs{G2}$, $\bs{G4}$.  Let
    $\bar{I}_{\theta}=I_{\theta} \backslash [\{ (1,0,0) , (0,1,0) \}
    \cup J_{|\la_1|} \cup J_{|\la_2|}]$.  Any trajectory
    $(y_1(s),y_2(s),\zeta(s))$ converging to $\bs{\mathcal{P}}$ is
    such that $y_1(s)$ has the following expansion if $l_s> \sigma^*$:
    \begin{equation}\label{expandy}
      \begin{aligned} & \ds y_1(s)= P_1^{+}+ \Psi(s) +a
        \eu^{\la_1 s}+ Q^1_{\theta}(s)+b \eu^{\la_2 s}+
        Q^2_{\theta}(s) +o( \eu^{-\theta s}) , \textrm{ where }\\
        & \ds Q_{1,\theta}(s)=\sum_{\chi \in
          J_{|\la_2|}} c^\chi \eu^{(\chi_1 \la_1+\chi_2 \la_2-\chi_3
          \gamma)s}, \textrm{ with $\chi=(\chi_1,\chi_2,\chi_3)$, and }
        \\  & \ds Q_{2,\theta}(s)=\sum_{\chi \in
          \bar{I}_{\theta}} c^{\chi} \eu^{(\chi_1 \la_1+\chi_2 \la_2-\chi_3
          \gamma)s}
      \end{aligned}
    \end{equation}
    as $s \to +\infty$, if we do not have resonances, otherwise we
    need to replace the resonant terms in $Q_{1,\theta}(s)$ according
    to \eqref{sost}.

    If $l_s=\sigma^*$ so that $\la_1=\la_2$ we have
    \begin{equation}\label{expandybis}
      \begin{array}{c}
        \ds y_1(s)= P_1^{+}+ \Psi(s)+a s\eu^{\la_1 s}+b \eu^{\la_1
          s}+ Q_{2,\theta}(s) +o( \eu^{-\theta s}) \,
      \end{array}
    \end{equation}
    as $s \to +\infty$, again if we do not have resonances, otherwise
    we need to replace the resonant terms in $ Q_{2,\theta}(s)$
    according to \eqref{sost}.
  \end{prop}

  \begin{remark}\label{importante}
    We emphasize that $Q_{1,\theta}(s)$ contains terms which are
    negligible with respect to $a \eu^{\la_1 s}$ while
    $Q_{2,\theta}(s)$ contains terms which are negligible with respect
    to $b \eu^{\la_2 s}$. Further if $|\la_1|<\gamma$ then $\Psi(s)$
    is identically null by definition.
  \end{remark}

  The proof is developed in the Appendix in a general
  framework, by showing a result on asymptotic expansions for ODEs,
  which seems to be new to the best of our knowledge.  In fact we
  borrow some of the ideas from \cite{DLL,YZ}.


  \begin{remark}\label{defa}
    Fix $\Q$ and $\tau \in \RR$; then $y_1(t,\tau,\Q;l_s)$ admits an
    expansion either of the form \eqref{expandy} or of the form
    \eqref{expandybis}.  All the coefficients in the expansions are
    determined by the choice of $a,b$, which are in fact smooth
    functions of $\Q$, i.e. $a=a(\Q)$, $b=b(\Q)$.

    In fact, all the coefficients in $\Psi(s)$ are determined when the
    non-linearity $g$ and $\tau$ are fixed; the coefficients in
    $Q_{1,\theta}$ are assigned (and can be determined) once $a$ is
    fixed, while $Q_{2,\theta}$ is assigned once $a$ and $b$ are
    assigned.
  \end{remark}

  \begin{remark}\label{datogliere}
    Fix $\Q$ and $\tau$, the coefficients $a=a(\Q)$, $b=b(\Q)$ may be
    evaluated through the method explained in \cite{DLL}. However from
    the previous discussion we have the following. Let $a_1,b_1, z_1$
    be such that $(\Q- \bs{P^+},\eu^{-\gamma \tau})=a_1
    v_1+b_1v_1+z_1v_z$.  Then $a=a_1+O(|\Q- \bs{P^+}|^2)$ and
    $b=b_1+O(|\Q- \bs{P^+}|^2)$.
  \end{remark}
  The proof of these two remarks is provided in
the Appendix.  For further details about these
  points 
  see \cite[Remarks~4.12, 4.16]{BF-arxiv}.

  Now, we translate Proposition \ref{sintetizzo} for the original
  equation \eqref{radsta}.
  \begin{lemma}\label{asinteqna}
    Assume $\bs{G1}$, $\bs{G2}$ with $l_u \ge l_s \ge \sigma^*$,
    $\bs{G3}$, $\bs{G4}$.  Consider a GS $U(r,\alpha)$ for $\alpha>0$,
    or the SGS $U(r,\infty)$ then there are continuous functions
    $\mathcal{A}\, :\, (0,+\infty] \to \RR$, $\mathcal{B}:(0,+\infty]
    \to \RR$, such that $\mathcal{A}$ is monotone decreasing, and if
    $l_s>\sigma^*$
    \begin{align}
      & \begin{aligned} \ds U(r,\alpha)= \frac{P_1^{+}}{r^{m}}+
        \frac{\Psi(\ln(r))}{r^m}+& \mathcal{A}(\alpha) r^{\la_1-m}+
        \frac{Q_{1,\theta}(\ln(r))}{r^m}\\ &+\mathcal{B}(\alpha)
        r^{\la_2-m} + \frac{Q_{2,\theta}(\ln(r))}{r^m} +o(
        r^{-\theta-m})
      \end{aligned} \label{expandU} \intertext{as $r \to +\infty$.  If
        $l_s=\sigma^*$ we have} &\begin{aligned} \ds U(r,\alpha)=
        \frac{P_1^{+}}{r^{m}}+
        \frac{\Psi(\ln(r))}{r^m}+\mathcal{A}(\alpha)\ln(r)
        &r^{\la_1-m} +\mathcal{B}(\alpha) r^{\la_2-m}\\ &+
        \frac{Q_{2,\theta}(\ln(r))}{r^m} +o( r^{-\theta-m}).
      \end{aligned} \label{expandUbis}
    \end{align}
  \end{lemma}

  \begin{remark} If we replace $\bs{G3}$ with the weaker assumption
    $\bs{A^-}$ in Lemma \ref{asinteqna}, then we still get the
    expansions in \eqref{expandU}, \eqref{expandUbis}, but we cannot
    ensure that $\mathcal{A}$ is monotone decreasing.
  \end{remark}

  \begin{proof}[ Proof of Lemma~\ref{asinteqna}]
    Fix $\tau \in \RR$; let $\bs{y}(s,\tau, \Q(\alpha);l_s)$ be the
    trajectory of \eqref{si.na} corresponding to $U(r,\alpha)$, so
    that $\Q(\alpha) \in W^u_{l_s}(\tau)$.  Then we can apply
    Proposition \ref{sintetizzo} to $y_1(s,\tau, \Q(\alpha);l_s)$ and
    we find the expansions \eqref{expandU}, \eqref{expandUbis}, where,
    according to Remark \ref{defa}, the coefficients $a,b$ are
    $a=a(\Q(\alpha))$ and $b=b(\Q(\alpha))$.  We set
    \begin{equation}\label{defAA}
      \mathcal{A}(\alpha)=a(\Q(\alpha)) \, , \qquad
      \mathcal{B}(\alpha)=b(\Q(\alpha)) \,.
    \end{equation}
    It follows that $\mathcal{A}:(0,+\infty) \to \RR$ and
    $\mathcal{B}:(0,+\infty) \to \RR$ are continuous functions.
    Finally if $\bs{G3}$ holds $U(r,\alpha_1)< U(r,\alpha_2)$ if
    $\alpha_1<\alpha_2$ for any $r>0$, and in particular for $r$
    large, so $\mathcal{A}(\alpha)$ is monotone increasing.
  \end{proof}

  \section{Main results: Stability and asymptotic
    stability} \label{sec:stability} Let us state
  Theorems~\ref{stabile} and \ref{asint.stabile} from which
  Theorems~\ref{thmD}, \ref{main-1}, \ref{main-2} follow directly.
  Let $r>0$, we denote by $[[r]]:=\{ k \in \mathbb{N} \mid k-1<r \le k
  \}$. We have the following results

\begin{thm}\label{stabile}
  Suppose $f$ is $C^k$ where $k=[[|\la_1|/\gamma]]$. Assume $\bs{K}$,
  $\bs{G1}$, $\bs{G2}$, $\bs{G3}$, $\bs{G4}$.  Then any radial GS
  $U(r,\alpha)$ of \eqref{parab} is stable with respect to the norm
  $\|\cdot\|_{m(l_s)+\lambda_1}$ if $l_s> \sigma^*$, and with respect
  to the norm $\vt{\,\cdot\, }_{m(l_s)+|\lambda_1|}$ if
  $l_s=\sigma^*$.
\end{thm}

\begin{thm}\label{asint.stabile}
  Assume the hypotheses of Theorem \ref{stabile}.  Then any radial GS
  $U(r,\alpha)$ of \eqref{parab} is weakly asymptotically stable with
  respect to the norm $\|\cdot\|_{m(l_s)+|\la_2|}$ if $l_s> \sigma^*$,
  and with respect to the norm $\| \, \cdot \,
  \|_{m(l_s)+|\lambda_1|}$ if $l_s=\sigma^*$.
\end{thm}
Let us recall that the stability of positive GS $U(|x|, \alpha)$ of
\eqref{parab} has been analyzed in a number of papers, (see
\cite{DLL,GNW1, GNW2, W}).  In \cite{GNW1}, when $f(u, |x|)=u^{q-1}$
and $q>\sigma^\ast$, the authors proved that the positive GS of
\eqref{parab} are stable in the norm $\|\, \cdot\, \|_{m+
  |\lambda_1|}$, and weakly asymptotically stable with respect to $\|
\, \cdot\, \|_{m+ |\lambda_2|}$. These results have been subsequently
extended in \cite{DLL} to functions $f(u, |x|)$ of the form
$k(|x|)r^{\delta}|u|^{q-1}$ where $K$ is a monotone decreasing
uniformly positive and bounded function. Here, we are able to prove
\emph{asymptotic stability} in place of \emph{weak asymptotic
  stability}. Further, we drop the assumption that $k$ is bounded:
This will allow us to consider potential giving rise to singular
solutions $U(r,\infty)$ having two different behaviors as $r \to 0$
(i.e. $U(r,\infty) \sim P^- r^{-m(l_u)}$) and as $r \to \infty$
(i.e. $U(r,\infty) \sim P^+ r^{-m(l_s)}$).

\subsection{Proof of Theorem \ref{stabile}.}
We first introduce some standard definition.
\begin{defin}
  We say that $\ov{\phi}$ is a super-solution of \eqref{laplace} if
  $\Delta \ov{\phi} + f(\ov{\phi},|x|) \le 0$; analogously
  $\und{\phi}$ is a sub-solution if $\Delta \und{\phi} +
  f(\und{\phi},|x|) \ge 0$.
\end{defin}
We refer to \cite{W} or to \cite[\S 3]{BF} for an extension of this
definition to weak and mild solutions.  Also, depending on a number of
very relevant factors (for instance, the type of domain and of the
boundary conditions, the regularity of the forcing term, etc... ) the
notion of weak solution for parabolic equations can change
considerably as described, e.g., in \cite{BH, Friedman, Gazzola,
  QS}. In particular, we mention that, a dynamical approach to study a
generalized parabolic equation on an unbounded strip-like domain is
given in \cite{Polat}: In this case a suitable definition of weak
solutions, on weighted Sobolev (and Bochner) spaces, is considered and
the author proved the existence of a global attractor.  Then, this
situation is further generalized in \cite{BC}.  \smallskip

Both Theorems \ref{stabile}, \ref{asint.stabile} depend strongly on
the following well known fact, proved in \cite[Theorem 2.4]{W}, see
also \cite[Theorem 3.10]{BF}.
\begin{lemma} \label{keyintrobis} Assume $\bs{G1}$, $\bs{G2}$, and let
  $U_1(r)$ and $U_2(r)$ be positive solutions of \eqref{radsta}
  respectively for $r \le R_1$ and for $r \ge R_2$, where $R_1>R_2$,
  and let $R \in (R_2,R_1) $ be such that $U_1(R)=U_2(R)$.  Consider
  \begin{equation*}
    \phi(x)=\left \{\begin{array}{ccl}
        U_1(r) & \textrm{if} & 0<|x| \le R, \\
        U_2(r) & \textrm{if} & |x| \ge R.
      \end{array} \right.
  \end{equation*}
  We have that
  \begin{itemize}
  \item If $U'_1(R)\ge U'_2(R)$, then $\phi(x)$ is a continuous weak
    super-solution of \eqref{laplace}.\\[-0.35 cm]
  \item If $U'_1(R)\le U'_2(R)$, then $\phi(x)$ is a continuous weak
    sub-solution of \eqref{laplace}.
  \end{itemize}
\end{lemma}
\begin{lemma} \label{keyintro} Assume $\bs{G1}$, $\bs{G2}$;
  \begin{itemize}
  \item [\textbf{(i)}] If the initial value $\phi$ in \eqref{data} is
    a continuous weak super-(sub-) solution of \eqref{laplace}, then
    the solution $u(t,x; \phi)$ of \eqref{parab}-\eqref{data} is
    non-increasing (non-decreasing) in $t$ as long as it exists, for
    any $x$;
    strictly if $\phi$ is not a solution.\\[-0.3 cm]
  \item[\textbf{(ii)}] If $\phi$ is radial, then $u(t,x; \phi)$ is
    radial in the $x$ variable for any $t>0$.\\[-0.3 cm]
  \end{itemize}
\end{lemma}
\noindent
To prove Theorem \ref{stabile} we adapt the main ideas developed in
\cite{GNW1,DLL,YZ}.

As a consequence of the proof of Proposition \ref{ord1-} we get the
following result which will be useful to prove the stability of the
solutions, and replaces a longer elliptic estimate performed in
\cite[Lemma 4.3]{DLL} and adapted in \cite{YZ,Dnew}.
We stress that in fact the proof in the critical case, considered in
\cite{Dnew}, suffers from a flaw.
\begin{lemma}\label{strict}
  Assume $\bs{K}$, $\bs{G1}$, $\bs{G2}$, $\bs{G3}$, $\bs{G4}$. Assume
  $\beta>\alpha$ then $\mathcal{A}(\beta)>\mathcal{A}(\alpha)$.
\end{lemma}
\begin{proof}
  Since $U(r,\beta)>U(r,\alpha)$ for any $r>0$ (see
  Lemma~\ref{ord1-}), we already know that $\mathcal{A}(\beta) \ge
  \mathcal{A}(\alpha)$, so we just need to prove that the inequality
  is strict.  Set
  $h(s)=[U(\eu^{s},\alpha_2)-U(\eu^{s},\alpha_1)]\eu^{(m(l_s)-\la_1)s}$,
  and, following the notation of Proposition \ref{ord1-},
  $Q(s)=\eu^{\la_1 s}$. Following the main line in the proof of
  Proposition \ref{ord1-} we see that $\dot{h}(s)=Z(s)/Q^2(s)$. In
  particular, from \eqref{integral} and \eqref{convex}, $\dot{h}(s)>0$
  for any $s \in \RR$.  Since $\lim_{s \to -\infty} h(s)=0$ we see
  that $h(s)>0$ for any $s \in \RR$, and $\lim_{s \to +\infty}
  h(s)>0$.

  If $l_s>\sigma^*$, then $\lim_{s \to +\infty}
  h(s)=\mathcal{A}(\beta)-\mathcal{A}(\alpha)>0$, and the proof is
  concluded.

  Assume now $l_s=\sigma^*$, and also assume by contradiction that
  $\mathcal{A}(\beta)= \mathcal{A}(\alpha)$. In this case we see that
  $\lim_{s \to +\infty} h(s)=\mathcal{B}(\beta)-\mathcal{B}(\alpha)
  \in (0, +\infty)$.  However, from \eqref{integral}, since
  $A=-2\la_1$, for any $\bar{s}\in \RR$ we find
  \begin{equation*}
    \begin{aligned}
      \dot{h}(\bar{s})= \int_{-\infty}^{\bar{s}} \eu^{A s}
      Q(s)[\partial_{y_1}g (P_1^+,+\infty;l_s)W(s)-D(s)]>0
    \end{aligned}
  \end{equation*}
  Therefore $\liminf_{s \to +\infty} \dot{h}(s)\ge \dot{h}(0) >0$,
  hence $\mathcal{B}(\beta)-\mathcal{B}(\alpha) =\lim_{s \to +\infty}
  h(s)=+\infty$, but this is a contradiction.  Hence
  $\mathcal{A}(\beta)> \mathcal{A}(\alpha)$.
\end{proof}

\begin{lemma}\label{vicinanza}
  Assume $\bs{K}$, $\bs{G1}$, $\bs{G2}$, $\bs{G3}$, $\bs{G4}$, and
  $l_u \ge l_s$.  If $l_s>\sigma^*$, then
  $\|U(r,\beta)-U(r,\alpha)\|_{m+ |\la_1|} \to 0$ as $\beta \to
  \alpha$, while if $l_s= \sigma^*$, then $ \vt{\,
    U(r,\beta)-U(r,\alpha) \,}_{m+ |\la_1|}\to 0$ as $\beta \to
  \alpha$.
\end{lemma}
\begin{proof}
  We develop the proof assuming $l_s>\sigma^*$, the case
  $l_s=\sigma^*$ is completely analogous.  It is well known that, for
  any fixed $R>0$ and any $\ep>0$, there is $\delta_1(\ep,\alpha,R)>0$
  such that
  \begin{equation}\label{azero}
    \sup \{|U(r,\beta)-U(r,\alpha)| \mid 0 \le r \le R \} < \ep
  \end{equation}
  whenever $|\beta-\alpha|<\delta_1$ (this is a continuous dependence
  on initial data argument for the singular equation \eqref{radsta}).
  Further from \eqref{expandU} we see that for $r$ large enough we
  have
  \begin{equation}\label{dentro}
    \begin{aligned}
      \left| \big(U(r,\beta)-U(r,\alpha)\big)\right|
      (1+r^{m-\lambda_1}) \cong
      |\mathcal{A}(\beta)-\mathcal{A}(\alpha)|+o(r^{|\la_2-\la_1|/2})
    \end{aligned}
  \end{equation}
  Thus, for any $\ep>0$ there exists $M(\ep)$ such that $
  |o(r^{|\la_2-\la_1|/2})|\leq C\ep $, when $r \ge M(\ep)$.  Further
  from Lemma \ref{asinteqna} we see that for any $\ep>0$ we can find
  $\delta_2(\ep,\alpha)>0$ such that $
  |\mathcal{A}(\beta)-\mathcal{A}(\alpha)| \le \ep$ if
  $|\beta-\alpha|<\delta_2$.  Therefore
  \begin{equation}\label{fuori}
    \begin{aligned}
      \left| \big(U(r,\beta)-U(r,\alpha)\big)\right|
      (1+r^{m-\lambda_1}) \le \ep \, , \quad \textrm{for $r \ge M$}
    \end{aligned}
  \end{equation}
  The proof then follows from \eqref{dentro}, \eqref{fuori}, choosing
  $M=R$ and $\delta(R,\alpha,\ep)= \min \{\delta_1, \delta_2 \}$.
\end{proof}


\begin{proof}[Proof of Theorem~\ref{stabile}]
  We give the proof just in the $l_s>\sigma^*$ case, in the
  $l_s=\sigma^*$ case is completely analogous.  Fix $\alpha>0$ and
  $\ep>0$ (small); let $\phi(x)$ be such that $\|U(|x|,
  \alpha)-\phi(x)\|_{m+|\la_1|} =\delta$, where $\delta>0$ will be
  chosen below.

  Let $|\eta|<\alpha$ and set
  \begin{equation}
    \begin{array}{c}
      z(r,\eta)= [U(r, \alpha+ \eta)-U(r,\alpha)](1+r^{m -\la_1})
    \end{array}
  \end{equation}
  Observe that $z(0,\eta)=\eta$ and $\lir z(r,\eta)=
  \mathcal{A}(\alpha+\eta)- \mathcal{A}(\alpha)$. So we can set
  \begin{equation}
    \underline{z}(\eta)= \min \{ |z(r,\eta)| \mid r>0 \} \,\, \textrm{
      and }\,\, \overline{z}(\eta)= \max \{ |z(r,\eta)| \mid r>0 \}.
  \end{equation}
  Moreover $z(r,\eta)$ is uniformly positive (respectively negative)
  for any $r>0$ if $\eta>0$ (resp. $\eta<0$), so
  $\underline{z}(\eta)>0$ if $\eta \ne 0$: This follows from Lemmas
  \ref{ord1-}, \ref{strict}.

  Finally, from Lemma ~\ref{vicinanza}, we know that $\lim_{\eta \to
    0} \underline{z}(\eta)= \lim_{\eta \to 0} \overline{z}(\eta)=0$.
  Then, for any $\ep>0$ we can find $d=d(\ep)>0$ such that
  $\overline{z}(-d) < \ep$, and $\overline{z}(d) < \ep$.  Set
  $\alpha_1=\alpha-d$, $\alpha_2=\alpha+d$, and choose $\delta= \min
  \{ \underline{z}(-d)\,, \underline{z}(d) \}$.  Then
  \begin{equation}
    \begin{aligned}
      U(|x|,\alpha_1) < \phi(x) < U(|x|,\alpha_2) & \, \\
      \|U(|x|,\alpha_i) -U(|x|,\alpha)\|_{m-\la_1} \le \ep & \textrm{
        for $i=1,2$}
    \end{aligned}
  \end{equation}
  Therefore, from the comparison principle (see, e.g.,
  \cite[Appendix]{Gazzola}), we have that
  \begin{equation}
    U(|x|,\alpha_1) < u(t,x;\phi) < U(|x|,\alpha_2), \textrm{ for
      any $t \ge 0$, $x \in \RR^n$},
  \end{equation}
  and the proof is concluded.
\end{proof}

%

\subsection{Weak asymptotic stability.}
To prove weak asymptotic stability we follow the outline of the proof
of \cite[Theorem 4.1]{GNW1} and adapted in \cite{DLL,YZ}.
\begin{prop}\label{4.1}
  Assume we are under the hypotheses of Theorem \ref{stabile} and
  consider the stationary problem \eqref{radsta}.  Then, for any
  radial GS $U(\cdot,d)$ of \eqref{radsta}, there is a sequence of
  radial strict super-solutions
  $\overline{U}^{(1)}(\cdot,e^1)>\overline{U}^{(2)}(\cdot,e^2)>\ldots>U(\cdot,d)$
  of \eqref{laplace} and a sequence of radial strict sub-solutions
  $\underline{U}^{(1)}(\cdot,c^1)<\underline{U}^{(2)}(\cdot,c^2)
  <\ldots<U(\cdot,d)$ such that $U(\cdot,d)$ is the only solution of
  \eqref{laplace} satisfying
  $\underline{U}^{(k)}(\cdot,c^k)<U(\cdot,d)<\overline{U}^{(k)}(\cdot,e^k)$,
  for every $k$.  Moreover
  \begin{equation}\label{limit}
    \lim_{k \to \infty}\underline{U}^{(k)}(\cdot,c^k)=U(\cdot,d)=
    \lim_{k \to \infty}\overline{U}^{(k)}(\cdot,e^k)
  \end{equation}
\end{prop}

  \begin{proof}
    Let $h:[0,+\infty) \to [0,1]$ be a monotone decreasing
    $C^{\infty}$ function such that $h(0)=1$ and $h(r)\equiv 0$ for $r
    \ge 1$.  Let
    $\mathcal{G}(y_1,s;l_s)=g(y_1,s;l_s)-g(y_1,+\infty;l_s)$ and
    observe that $\mathcal{G}(y_1,s;l_s) \ge 0$ and it is decreasing
    in $s$
    for any $y_1,s$.\\
    \noindent \textbf{--Assume first $\mathcal{G}(y_1,s) \not\equiv
      0$}, i.e. consider the generic case, and denote by
    \begin{equation*}
      \begin{array}{l}
        \ov{g}^{(k)}(y_1,s)= g(y_1,s;l_s)+\frac{h(\eu^s)}{2k}
        \mathcal{G}(y_1,s;l_s) \\ \und{g}^{(k)}(y_1,s)=
        g(y_1,s;l_s)-\frac{h(\eu^s)}{2k} \mathcal{G}(y_1,s;l_s)
      \end{array}
    \end{equation*}
    and let $\ov{f}^{(k)}$, $\und{f}^{(k)}$ be the corresponding
    functions obtained via \eqref{transf1}.  Notice that by
    construction $\ov{g}^{(k)}(y_1,s)$, and $\und{g}^{(k)}(y_1,s)$ are
    both decreasing in $s$ for any $k \ge 1$; hence $\ov{f}^{(k)} \ge
    f \ge \und{f}^{(k)}$ satisfy $\bs{G1}$, $\bs{G2}$, $\bs{G3}$,
    $\bs{G4}$, $\bs{K}$ so that Lemma~\ref{ord0}, and
    Proposition~\ref{ord1-} hold true.  In particular all the regular
    solutions of the respective problem \eqref{radsta}, say
    $\ov{U}^{(k)}(r,\alpha)$, $U(r,\alpha)$,
    $\und{U}^{(k)}(r,\alpha)$, are GSs. Further the corresponding
    trajectories of \eqref{si.na}, say $\bs{\ov{y}^{(k)}}(s,\alpha)$,
    $\bs{y}(s,\alpha)$, $\bs{\und{y}^{(k)}}(s,\alpha)$ are monotone
    increasing in their first component and converge to $\bs{P^+}$,
    and have the asymptotic expansion as described in Proposition
    \ref{sintetizzo}.  More precisely they both have either the
    expansion \eqref{expandU} or \eqref{expandUbis}, where the
    function $\Psi(\ln(r))$ coincide for $r \ge 1$, while the
    coefficients $a=\ov{\mathcal{A}}^{(k)}(\alpha)$,
    $a=\und{\mathcal{A}}^{(k)}(\alpha)$ and
    $b=\ov{\mathcal{B}}^{(k)}(\alpha)$,
    $b=\und{\mathcal{B}}^k(\alpha)$ are different, see
    Lemma~\ref{asinteqna}.  Further by construction,
    $\ov{U}^{(k)}(r,\alpha)$, $\und{U}^{(k)}(r,\alpha)$ are
    respectively super and sub-solutions for the original problem
    \eqref{radsta}.

    We divide our argument in several steps.  \smallskip

    \noindent \textbf{--Step 1.} $\;$\emph{If there is $R>0$ such that
      $U(R,d)=\ov{U}^{(k)}(R,c)$ (respectively
      $U(R,d)=\und{U}^{(k)}(R,e)$), then $U(r,d) \ge
      \ov{U}^{(k)}(r,c)$ (respectively $U(r,d)\le \und{U}^{(k)}(r,e)$)
      for
      any $r \ge R$.}\\
    Let $\tau(\xi): (0,2) \to \RR$ be the inverse of the function
    $\xi(\tau)$ defined in \eqref{manistella}. We consider
    \begin{equation}\label{si.nastella}
      \begin{pmatrix}
        \dot{y}_{1} \\ \dot{y}_{2} \\ \dot{z}
      \end{pmatrix}
      = \begin{pmatrix} 0 & 1 &0 \\ B(l_*) & -A(l_*) & 0 \\ 0 & 0 & \varpi
      \end{pmatrix}
      \begin{pmatrix}
        y_1 \\ y_2 \\ \xi
      \end{pmatrix}-
      \begin{pmatrix}
        0 \\ g(y_1,\tau(\xi);l_*)\\ 0
      \end{pmatrix},
    \end{equation}
    where $A(l_*)$, $B(l_*)$ coincide with $A(l_u)$, $B(l_u)$ for
    $s\le 0$ and with $A(l_s)$, $B(l_s)$ for $s\ge 0$, and similarly
    $g(y_1,\tau(\xi);l_*)$ equals $g(y_1,\frac{\ln(\xi)}{\varpi};l_u)$
    for $\xi \le 1$ (i.e. $s \le 0$) and
    $g(y_1,\frac{\ln(2-\xi)}{\varpi};l_s)$ for $\xi \ge 1$ (i.e. $s
    \ge 0$). Notice that \eqref{si.nastella} coincides with
    \eqref{si.naa} when $\xi \le 1$ (i.e. $s \le 0$) and it is
    equivalent to \eqref{si.naas} when $\xi \ge 1$ (i.e. $s \le 0$ and
    $\zeta \le 1$, it differs from \eqref{si.naas} just in the fact
    that $\xi=2-\zeta$).  Further we recall that the unstable manifold
    $\bs{W^u}(l_*)$ defined in \eqref{manistella} has dimension $2$
    and connects the $\xi$-axis and the graph of $\bs{y^{*}}(s,l_*)$];
    further it is a graph on the $y_2=0$ plane, see
    Remark~\ref{ord3-bis}.  Moreover $\bs{W^u}(l_*)$ splits the set
    \begin{equation*}
      E=\{(y_1,y_2,\xi) \mid 0<y_1<y_1^u(\tau(\xi),l_*) \, , \;
      0<\xi< 2 \}
    \end{equation*}
    in $2$ open components, say $E^+$ and $E^-$ (the one with larger
    and smaller $y_2$).

    By construction the flow of the modified system
    \eqref{si.nastella} where $g$ is replaced respectively by
    $\ov{g}^{k}$ and by $\und{g}^{k}$ on $\bs{W^u}(l_*)$ points
    towards $E^-$ and $E^+$ respectively for $s \le 0$, and it is
    tangent to $E^0$ for $s \ge 0$.  So the corresponding manifolds
    $\bs{\ov{W}^{u,(k)}}(l_*)$ and $\bs{\und{W}^{u,(k)}}(l_*)$ lie
    respectively in $E^-$ and $E^+$.

    Now assume $U(R,d)=\ov{U}^{(k)}(R,c)$ and consider the
    corresponding trajectories $\bs{y}(s;l_*)$, and $\bs{\ov{y}^{(k)}
    }(s;l_*)$: Then $y_{1}(\ln(R); l_*)=\ov{y}^{(k)}_{1}(\ln(R);l_*)$
    and $y_{2}(\ln(R); l_*) \ge \ov{y}^{(k)}_2(\ln(R);l_*)$.  Hence
    $y_{1}(s; l_*)\ge \ov{y}^{(k)}_{1}(s;l_*)$ for $s$ in a right
    neighborhood of $\ln(R)$.  Then the claim in Step 1 concerning
    $\ov{U}^{(k)}(r,c)$ follows. The claim concerning
    $\und{U}^{(k)}(r,e)$ is analogous.

    We continue the discussion for later purposes.  We know that
    $y_{2}(\ln(R); l_*) \ge \ov{y}^{(k)}_2(\ln(R);l_*)$, assume first
    $y_{2}(\ln(R); l_*) > \ov{y}^{(k)}_2(\ln(R);l_*)$. Then $y_{1}(s;
    l_*) > \ov{y}^{(k)}_1(s;l_*)$ for $s$ in a right neighborhood of
    $\ln(R)$.

    Assume now $y_{2}(\ln(R); l_*) = \ov{y}^{(k)}_2(\ln(R);l_*)$: Then
    $R \ge 1$. In fact assume for contradiction that $0<R<1$, then
    $\bs{\ov{y} }( \ln(R);l_*)=\bs{Q}=\bs{\ov{y}^{(k)} }( \ln(R);l_*)$
    is such that $(\Q,\xi(\ln(R))\in E^0$, but from \eqref{si.na} we
    get $\dot{\ov{y}}_{2}(\ln(R); l_*) <
    \dot{y}_2^{(k)}(\ln(R);l_*)$. Hence $\bs{\ov{y}^{(k)} }( r;l_*)$
    crosses transversally $E^0$ at $s=\ln(R)$, going from $E^+$ to
    $E^-$, in particular it is in $E^+$ when $s$ is in a sufficiently
    small left neighborhood of $\ln(R)$.  But
    $(\bs{\ov{y}^{(k)}}(s;l_*), \xi(s)) \in
    \bs{\ov{W}^{u,(k)}}(l_*)\subset E^-$, and this is a contradiction,
    so $R>1$.

    Observe that if $R \ge 1$ then $\bs{\ov{y}^{(k)}}(s;l_*)$ and
    $\bs{y}(s;l_*)$ are solutions of the same equation \eqref{si.na}
    for $s \ge 0$ which coincide for $s=\ln(R)$, so they coincide for
    $s \ge 0$.

    In fact we have already proved the following, i.e. \smallskip

    \noindent \textbf{--Step 2.}  For any $0<r<1$ we have that
    \begin{equation}\label{ineq}
      \ov{U}^{(k)}(r,d) < U(r,d) < \und{U}^{(k)}(r,d)
    \end{equation}
    \emph{and either \eqref{ineq} holds for any $r>0$ or the functions
      coincide for any $r \ge 1$.  Moreover $\ov{\mathcal{A}}^{(k)}(d)
      \le \mathcal{A}(d) \le \und{\mathcal{A}}^{(k)}(d)$.}  \smallskip

    \noindent \textbf{--Step 3.} \emph{Fix $d$ and the corresponding
      coefficient $\mathcal{A}(d)$.  It is possible to choose $c^k\le
      d \le e^k$ so that $\und{\mathcal{A}}^{(k)}(e^k)=
      \ov{\mathcal{A}}^{(k)}(c^k)=\mathcal{A}(d)$. Then, from
      \emph{Step 1} it follows that} $U(r,d)$ is the unique solution
    of the original equation \eqref{radsta} such that
    \begin{equation}\label{eq:step3}
      \und{U}^{(k)}(r,c^k) \le U(r,d) \le \ov{U}^{(k)} (r,e^k) , \,\,
      \textrm{for any $r \ge 0$}
    \end{equation}

    Fix $\tau >0$ and $0<c<d<e$; let $\bs{y}(s,\tau,\bs{P}; l_s)$,
    $\bs{y}(s,\tau,\bs{Q}; l_s)$, $\bs{y}(s,\tau,\bs{R}; l_s)$ be the
    trajectories of \eqref{si.na} corresponding to the solutions
    $U(r,c)$, $U(r,d)$, $U(r,e)$ of \eqref{radsta}.  It follows that
    $\bs{P}, \Q, \bs{R}$ are points in $W^u(\tau,l_s)$ and $\bs{P}$,
    $\bs{R}$ are respectively the closest to and the farthest from the
    origin.  Let us consider the lines $\ell^l, \ell^r$ parallel to
    the $y_2$-axis and passing through $\bs{P}$ and $\bs{R}$
    respectively: We denote by $\bs{\ov{P}^{(k)}}$ and
    $\bs{\ov{R}^{(k)}}$, the intersections of
    $\ov{W}^{u,(k)}(\tau,l_s)$ respectively with $\ell^l$ and with
    $\ell^r$.  Using continuous dependence on initial data of ODE we
    see that $\bs{\ov{P}^{(k)}} \to \bs{P}$ and $\bs{\ov{R}^{(k)}} \to
    \bs{R}$ as $k \to \infty$.  Since $a(\Q)$ is continuous, see
    Remark \ref{defa} and \eqref{defAA}, we see that
    $a(\bs{\ov{P}^{(k)}}) \to
    a(\bs{P})=\mathcal{A}(c)<\mathcal{A}(d)$, while
    $a(\bs{\ov{R}^{(k)}}) \to
    a(\bs{R})=\mathcal{A}(e)>\mathcal{A}(d)$. Therefore we can choose
    $N$ large enough so that
    $a(\bs{\ov{P}^{(k)}})<\mathcal{A}(d)<a(\bs{\ov{R}^{(k)}})$ for any
    $k \ge N$. Hence we can find $\bs{\ov{Q}^{(k)}}\in
    \ov{W}^{u,(k)}(\tau,l_s)$ between $\bs{\ov{P}^{(k)}}$ and
    $\bs{\ov{R}^{(k)}}$ such that
    $a(\bs{\ov{Q}^{(k)}})=\mathcal{A}(d)$.
    Correspondingly we find $e^k$ such that
    $\ov{\mathcal{A}}^{(k)}(e^k)=a(\bs{\ov{Q}^{(k)}})=\mathcal{A}(d)$.
    Note that in view of \emph{Step 2} we have $e^k \ge d$.  The proof
    for $\und{\mathcal{A}}^{(k)}(c^k)$ is analogous.  \smallskip

    \noindent \textbf{--Step 4.} \emph{Formula \eqref{limit} and the
      following Remark hold true.}

        \begin{remark}\label{estimateL}
          $\ov{\mathcal{B}}^{(k)}(e^k)$ and
          $\und{\mathcal{B}}^{(k)}(c^k)$ are respectively strictly
          decreasing and increasing in $k$ and they both converge to
          $\mathcal{B}(d)$.

        \end{remark} 
        \begin{proof}
          To prove \eqref{limit} it is enough to observe that, by
          construction, the functions $\und{U}^k(r,c^k)$ and
          $\ov{U}^k(r,e^k)$ are bounded and monotonically respectively
          increasing and decreasing in $k$. Then, from standard
          elliptic estimates and \emph{Step 3} we see that the limit
          of both is the solution $U(r,d)$ of the original problem
          \eqref{radsta}.

          Now, we turn to consider Remark \ref{estimateL}.  Let us
          recall that, by construction, the following relation holds
          true, i.e.:
          $\underline{\mathcal{A}}^{(k)}(e^k)=\mathcal{\mathcal{A}}(d)=
          \overline{\mathcal{A}}^{(k)}(c^k)$ (see \emph{Step
            3}). Further we also infer that
          $\ov{\mathcal{B}}^{(k)}(e^k)$ and
          $\und{\mathcal{B}}^{(k)}(c^k)$ are respectively decreasing
          and increasing and converge to $\mathcal{B}(d)$.  As next
          step, we show that $\und{\mathcal{B}}^{(k)}(c^k)<
          \und{\mathcal{B}}^{(k-1)}(c^{k-1})<\mathcal{B}(d)<
          \ov{\mathcal{B}}^{(k-1)}(e^{k-1})<\ov{\mathcal{B}}^{(k)}(e^k)$.
          As usual we just prove the last inequality, the others being
          analogous.  Let $\ov{u}^j(x)$ be the radial function defined
          by $\ov{u}^j(x)=\ov{U}^j(|x|,e^j)$. Observe that $\Delta
          [\ov{u}^k(x))-\ov{u}^{k-1}(x)] \le 0$, hence from standard
          arguments (see \cite[Theorem 3.8]{Ni}),
          we see that there is $C>0$ such that
          $\ov{U}^k(r,e^k)-\ov{U}^{k-1}(r,e^{k-1})> C r^{-(n-2)}$.
          Assume
          $\ov{\mathcal{B}}^{(k)}(e^k)=\ov{\mathcal{B}}^{(k-1)}(e^{k-1})$
          for contradiction.  Since
          $\ov{\mathcal{A}}^{(k)}(e^k)=\mathcal{\mathcal{A}}(d)$ for
          any $k$, from the construction in Lemma \ref{asinteqna} it
          follows that $\bs{\ov{y}^{(k)}}(s, e^k;l_s) \equiv
          \bs{\ov{y}^{(k-1)}}(s, e^{k-1};l_s)$ for any $s \ge
          0$, 
          i.e.  $\ov{U}^k(r,e^k)= \ov{U}^{k-1}(r,e^{k-1})$ for $r \ge
          1$, but this is a contradiction and the Lemma is proved.
        \end{proof}

        From Remark \ref{estimateL} we see that the inequalities in
        \eqref{eq:step3} are strict for $r$ large.  Then, from
        \emph{Step 1} we conclude.\smallskip

        \textbf{--Assume now $G(y_1,s) \equiv 0$}, this is the case,
        e.g., when $f(u,r)= c u|u|^{q-2}$.  Following \cite{GNW1} we
        denote by $\underline{f}^{(k)}(u,r):= [1- \mu h(r)/k]f(u,r)$
        and $\overline{f}^{(k)}(u,r):= [1+ \mu h(r)/k]f(u,r)$, for $k
        \in \mathbb{N}$ and where $\mu>0$ is chosen small enough so
        that $\und{f}^{(1)}(u,r)$ satisfies $\bs{A^-}$; then it is
        easy to check that $\und{f}^{(k)}(u,r)$ and
        $\ov{f}^{(k)}(u,r)$ satisfy $\bs{A^-}$ for any $k \in
        \mathbb{N}$. So Proposition~\ref{super} holds true, and in all
        the $3$ cases all the regular solutions of \eqref{radsta},
        denoted respectively by $\und{U}^{(k)}(r,\alpha)$,
        $U(r,\alpha)$, $\ov{U}^{(k)}(r,\alpha)$, are GSs, but a a
        priori they might not be ordered.  However repeating the
        argument of \emph{Step 1} in \cite[Theorem 4.1]{GNW1}, it it
        easy to prove that
        \begin{equation}\label{pre-step-2}
          \ov{U}^{(k)}(r,\alpha) \le U(r,\alpha) \le
          \und{U}^{(k)}(r,\alpha)
        \end{equation}
        for any $r >0$ and any $\alpha > 0$.  The proof might be
        concluded arguing as in \cite[Theorem 4.1]{GNW1}.  However
        notice that we can also repeat the argument at the end of
        \emph{Step 1} of this proof to get \eqref{ineq} for any $r>0$,
        and then carry on through \emph{Step 2,3,4} of this proof and
        conclude also in this case, with no further changes.
      \end{proof}
      From the previous discussion we easily find the following
      result.

      \subsection{Proof of the weak asymptotic stability}
%
      Now we consider $d>0$ fixed, and we use the shorthand notation
      $\overline{U}^{(1)}(r,e^1)= \overline{U}(r)$,
      $\underline{U}^{(1)}(r,c^1)= \underline{U}(r)$,
      $\overline{u}(t,x)=u(t,x; \overline{U}(|x|))$,
      $\underline{u}(t,x)=u(t,x; \underline{U}(|x|))$.
      \begin{lemma}\label{forteasintotico}
       Assume that we are in the
        hypotheses of Theorem \ref{stabile}; Then $\overline{u}(t,x)
        \searrow U(|x|,d)$ and $\underline{u}(t,x) \nearrow U(|x|,d)$
        as $t \to +\infty$, with the norm $\| \cdot \|_{ l}$, for any
        $0 \le l<m+|\la_2|$.
      \end{lemma}
      Notice that if $l_s=\sigma^*$ then $\| \cdot \|_{ m+|\la_1|}=\|
      \cdot \|_{ m+|\la_2|}$.
      \begin{proof}
        Let us set $B:=\lim_{|x| \to
          +\infty}[\overline{U}(|x|)-\underline{U}(|x|)]|x|^{m+|\la_2|}$
        and notice that $B>0$ is finite, see Proposition \ref{4.1} and
        Remark \ref{estimateL}.  Fix $0 \le l<m+|\la_2|$ and observe
        that for any $\ep>0$ we can find $\rho>0$ such that
        \begin{equation}\label{Due}
          [\overline{U}(|x|)-\underline{U}(|x|)]|x|^{l}<2 B |x|^{l-m+|\la_2|}< \ep/2
        \end{equation}
        for $\|x \| \ge \rho$.

        Since $\overline{U}(|x|)$ and $\underline{U}(|x|)$ are
        respectively a radial super and sub-solution of
        \eqref{laplace}, then $\overline{u}(t,x)$ and
        $\underline{u}(t,x)$ are respectively radially symmetric super
        and sub-solution of \eqref{parab}.  Further they are resp.
        monotone decreasing and increasing in $t$, so they converge to
        a radial solution of \eqref{laplace}, see Lemma
        \ref{keyintro}.  From Lemma \ref{4.1} we know that $U(r,d)$ is
        the unique solution of \eqref{radsta} between
        $\overline{U}(r)$ and $\underline{U}(r)$, so
        $\overline{u}(t,x)$ and $\underline{u}(t,x)$ converge
        monotonically to $U(|x|,d)$ as $t \to +\infty$, for any fixed
        $x \in \RR^n$.  Then, from the equiboundedness of the
        functions involved and of their derivatives we see that the
        convergence is uniform in any ball of radius $R>0$
        fixed. Hence setting $R=\rho>0$, for any $\ep>0$ we find
        $T(\ep)>0$ such that
        \begin{equation}\label{Uno}
          [\overline{u}(x,t)-\underline{u}(x,t)] |x|^l \le \ep/2
        \end{equation}
        for any $|x| \le \rho$. Further from \eqref{Due} and the
        comparison principle we easily find
        \begin{equation}\label{Tre}
          [\overline{u}(x,t)-\underline{u}(x,t)] |x|^l
          \le   [\overline{U}(|x|)-\underline{U}(|x|)]|x|^{l} \le \ep/2
        \end{equation}
        for $|x| \ge \rho$. Hence the Lemma follows from \eqref{Uno}
        and \eqref{Tre}.
      \end{proof}

    \begin{proof}[Proof of Theorem~\ref{asint.stabile}]
      Assume for definiteness $l_s>\sigma^*$, the case $l_s \ge
      \sigma^*$ being analogous.  Fix $d>0$ and denote by
      \begin{equation}\label{www}
        \begin{array}{cc}
          \overline{W}(r,d)= [\overline{U}(r)-U(r,d)] (1+r^{m+|\la_2|})
          & \overline{\delta}= \inf_{r>0}\overline{W}(r,d)
          \\ \underline{W}(r,d)= [U(r,d)-\underline{U}(r)]
          (1+r^{m+|\la_2|}) & \underline{\delta}=
          \inf_{r>0}\underline{W}(r,d)
        \end{array}
      \end{equation}
      Observe that $\overline{W}(r,d)$, $\underline{W}(r,d)$ are both
      positive for any $r>0$, see Proposition \ref{4.1}.  Further
      $\overline{W}(0,d)= e^1-d>0$, $\underline{W}(0,d)= d-c^1>0$,
      $\lir \overline{W}(r,d)= \overline{B}^{(1)}(e^1)-B(d)>0$, $\lir
      \underline{W}(r,d)=B(d)-\underline{B}^{(1)}(c^1)>0$, see Remark
      \ref{estimateL}. It follows that $\delta= \min \{
      \overline{\delta}, \underline{\delta} \}>0$.

      Now let us consider $\phi$ such that $\|\phi-U(| \cdot |,
      d)\|_{m+|\la_2|} < \delta$: by construction we have
      $\underline{U}(|x|) \le \phi(x) \le \overline{U}(|x|)$, for any
      $x \in \RR^n$. Therefore
    $$ \underline{u}(t,x) \le u(t,x ; \phi) \le \overline{u}(t,x)$$
    for any $t>0$ and any $x \in \RR^n$. So from Lemma
    \ref{forteasintotico} we easily conclude.
  \end{proof}

\appendix
 \section{Proof of Proposition~\ref{sintetizzo}}
  \label{sec:appendix}
   In what follows  we develop a constructive argument to prove
   the  asymptotic expansion results of subsection~\ref{subsec:asympotic}. This
   result, to the best of our knowledge, seems to be new  and it is
   of independent interest for the ODEs theory.
 Further, we borrow some of the ideas from
  \cite{DLL,YZ}.
  The  purpose is to approximate a generic solution converging
  to a critical point of a non-linear system, by a recursive sequence
  of solutions of approximating linear non-homogeneous systems.

  Let us consider an equation of the form
  \begin{equation}\label{sistemacompleto}
    \dot{x}=f(x)=L x+N(x)
  \end{equation}
  where $x \in \RR^n$,  $N(x)$ is at least  $C^2$ and
 such that $N(0)=N_x(0)=0$. More regularity will be required in the
 second part of the proof.

  We denote by $\la_i$ the eigenvalues of $L$ and by  $m$ is the number of  eigenvalues with distinct real parts. We set
  $\wedge_i:=-\textrm{Re}(\la_i)$, and we assume for definiteness
  $\wedge_i< \wedge_{i+1}$ for any $1 \le i \le m-1$. We also assume
  $\wedge_1>0$, hence \eqref{sistemacompleto} is exponentially stable.
  Further let $l_i$ be the number of eigenvalues (counted with multiplicity) such that $-\textrm{Re}(\la_j)=\wedge_i$, so that $\sum_{i=1}^{m} l_i=n$.

 We can assume
without loss of generality  that $L$ is block diagonalized,
and that each block, $L_1,
\ldots, L_m$ is in Jordan form.
We denote by $P_i$
the matrix which is the identity in the $i$-th block and the null matrix in the other blocks,
so that $L P_i=   P_i L=L_i$. Then we set $\PP_i=\sum_{j=1}^i P_j$:
So the matrices $P_i$ and $\PP_i$ are projections
on eigenspaces.

Note that any solution $\ell(t)$ of the linear equation   $\dot{y}= L y$  takes the form
\begin{equation}\label{blocks}
  \begin{aligned}
    &\ell(t)= (\eu^{L_1 t}, \eu^{L_2 t}, \ldots ,
    \eu^{L_m t}) \ell(0) = \sum_{i=1}^m\ell_i(t)\\[-0.15 cm]
    &\textrm{where} \qquad
    \ell_i(t)= (\underbrace{0,\ldots, 0}_{l_1+ \ldots + l_{i-1}},
    \eu^{L_i t}, \underbrace{0 \ldots, 0)}_{l_{i+1}+ \ldots + l_{m}}
    \ell(0)
  \end{aligned}
\end{equation}
In the whole appendix the notation $f(t)=o(g(t))$ and $f(t)=O(h(t))$
  means respectively that $\|f(t)\|/\|g(t)\| \to 0$ and $\|
  f(t)\|/\|h(t)\| \|$ remains bounded as $t \to +\infty$. Further
  $\ep$ is a positive constant which is taken as small as needed, and
  it
 may change from line-to-line.
  \begin{remark}\label{rem1}
    Notice that, for any $\ep>0$ we have
     $\eu^{-\wedge_i t} \le \|\eu^{L_i
  t}\| \le \eu^{(-\wedge_i+\ep) t}$
 whence $t>0$, and $i=1, \ldots,
m$. In fact  there is $c>0$ such that
   $ \eu^{-\wedge_i t} \le \|\eu^{L_it}\|   \le
  c [1+t^{l_i-1}]\eu^{-\wedge_i t}$ for any $t>0$, and any $i=1, \ldots, m$.
\end{remark}
  We start from the following technical Lemma.
  \begin{lemma}\label{stimazero}
  Let $y(t)=O(\eu^{-k t})$ where $k> \wedge_i$, then the integral
  $$I(t)=\int_t^{+\infty} \eu^{L(t-s)} \PP_i y(s) ds $$
  is well defined and $I(t)=O(\eu^{-k t})$.
  \end{lemma}
  \begin{proof}
 First notice that, for any $j=1, \ldots, i$, and any sufficiently small $\ep>0$, we have
\begin{equation*}
\begin{split}
\left \|\int_{t}^{+\infty} \eu^{L(t-s)}  P_j  y(s) ds \right \|  \le  \int_{t}^{+\infty} \eu^{(\wedge_j+\ep) (s-t)} C \eu^{- k s} ds
   =  \frac{C}{k-\wedge_j-\ep} \eu^{-k t}
\end{split}
\end{equation*}
Then observe that the first $i$ blocks of $I(t)$ satisfy the previous
estimate,
while the last $m-i$ ones are null. So
the Lemma immediately follows.
  \end{proof}
Now we recall the
  following standard result: We sketch the proof since it gives the
\textbf{Step 0} of our approximating procedure.

  \begin{lemma}\label{linn}\cite[\S 13-Theorem 4.5]{CoLe}
 Let  $f$ be $C^2$, and let $x(t)$ be a solution of
 \eqref{sistemacompleto}
 such that $x(t) \to 0$ as $t \to +\infty$.  Then,  there is $\bar{\ell}$ such that the
    solution $\ell(t):= \eu^{L t} \bar{\ell}$ of the linear equation
    $\dot{y}= L y$,  satisfies
    \begin{equation}\label{est.lin}
      \|x(t)-\ell(t)\| \eu^{ 2(\wedge_1 -\ep)t } \to 0 \, , \quad
      \textrm{as $t \to +\infty$,}
    \end{equation}
     for any $\ep>0$. Further
       \begin{equation}\label{N0}
    \bar{\ell}=x(0)+\bar{\mathfrak{N}}_1, \; \textrm{where} \; \; \bar{\mathfrak{N}}_1= \int_0^{+\infty} \eu^{-Ls} \PP_1 N(x(s))ds .
  \end{equation}
  \end{lemma}
\begin{proof}[Sketch of the proof]
  From \cite[\S 13]{CoLe} we already know that any solution $x(t)$ of
  \eqref{sistemacompleto} satisfies $x(t) \eu^{(\wedge_1-\ep)t} \to 0$
  as $t \to +\infty$, for any $\ep>0$.  Therefore the integral in \eqref{N0} defining $\bar{\mathfrak{N}}_1$ is convergent and
  $\bar{\ell}$ is well defined, see Lemma \ref{stimazero}. Then observe that the solution
  $x(t)$ of \eqref{sistemacompleto} can be rewritten as follows
  \begin{equation}\label{w.nonlin}
    x(t)= \eu^{L t} x(0) + \int_0^t \eu^{L(t-s)} N(x(s)) ds
    \, .
  \end{equation}
  Hence we get
    \begin{equation}\label{w.corretto}
    x(t)-\ell(t)=  \int_0^t \eu^{L(t-s)}(\Id-\PP_1) N(x(s)) ds -\int_t^{+\infty} \eu^{L(t-s)} \PP_1 N(x(s)) ds
    \, .
  \end{equation}
  Since $N(x(t))=o(\eu^{-2(\wedge_1 - \epsilon)t})$, using also Lemma \ref{stimazero} we get
\eqref{est.lin} and conclude the proof.
\end{proof}

\smallskip
\begin{defin}\label{d:defkappa}
Let us set $k_1=1$ and let $k_2 \in \mathbb{N}$, $k_{2} \ge 0$ such that $(k_{2}+k_1)\wedge_1
 \le \wedge_2 < (k_{2}+k_1+1) \wedge_1$, and,
 for later purposes, denote by $k_{i}$ the unique integer such that
 \begin{equation}\label{defkappa}
\wedge_1 \big( \sum_{j=1}^{i} k_{j}\big)
 \le \wedge_i < \wedge_1 \big(1+ \sum_{j=1}^{i} k_{j}\big) \quad \textrm{for $i=1, \ldots, m$.}
\end{equation}
\end{defin}

\noindent
The step $1$ in our approximating scheme will be to prove the following Lemma.
\begin{lemma}\label{l:step.a}
 Set $\bar{\ell}_1:= P_1 (x(0))+ \bar{\mathfrak{N}}_1$ and
 $a_{1}^{0 }(t):=\ell_1(t)=\eu^{L t} \bar{\ell}_1$.
 Then, we can expand $x(t)$ as follows:
 \begin{equation}\label{expand.a.fin}
 \begin{split}
 &  x(t)= A^{k_2}_1(t)+  R^{k_2}_1(t) \, , \quad \textrm{where }  \; A^{k_2}_1(t)= \sum_{j=0}^{k_2} a^j_1(t) \\
 & a^j_1(t)=o(\eu^{-[(1+j) \wedge_1-\ep]t}) \, , \quad \textrm{and }  \; R^{k_2}_1(t)=o(\eu^{-( \wedge_2-\ep)t})
 \end{split}
 \end{equation}
\end{lemma}
To help the reader with the notation we emphasize that
the apex indicates the step of the iteration, while the subscript indicates
the eigenvalue we are dealing with.
\begin{proof}
Let us start from  $A_{1}^{0 }(t)=a_{1}^{0 }(t):=\ell_1(t)$,
$R_{1}^{0 }(t):=x(t)-\ell_1(t)$. Repeating the computation of Lemma \ref{linn} (with $\bar{\ell}_1$ replacing $\bar{\ell}$)
we see that
\begin{equation}\label{def.R1}
a_{1}^{0 }(t)=o(\eu^{-(\wedge_1-\ep)t }), \quad      R_{1}^{0 }(t)=o(\eu^{-(2\wedge_1-\ep)t }+ \eu^{-(\wedge_2-\ep)t }).
\end{equation}
If $k_2=0$ we have  $R_{1}^{0 }(t)=o( \eu^{-(\wedge_2-\ep)t })$ and
Lemma \ref{l:step.a} is proved.
Otherwise we go through the following steps:
\smallskip


\noindent \textbf{--Step 1. } Since $N(a_{1}^{0 }(t))=o(\eu^{(-2\wedge_1+\ep)t
})$ as $t \to +\infty$, using Lemma \ref{stimazero},
we can define
\begin{equation}\label{def.a2}
  \begin{split}
  M_{1}^{1 }(s):= &  N(A_{1}^{0 }(s))  \, , \qquad \qquad \textrm{and} \\
 a_{1}^{1 }(t):=& \int_{0}^{t}  \eu^{L(t-s)}(\Id-\PP_1) M_{1}^{1 }(s) ds  - \int_{t}^{+\infty} \eu^{L(t-s)}  \PP_1 M_{1}^{1 }(s) ds
  \end{split}
\end{equation}
By construction
\begin{equation}\label{asymp.a2}
  a_{1}^{1 }(t)=o(\eu^{-2(\wedge_1-\ep)t}) \,\,  .
\end{equation}
Denote by $A_{1}^{1 }(t)=a_{1}^{0 }(t)+a_{1}^{1 }(t)$, and by $R_{1}^{1 }(t)=x(t)-A_{1}^{1 }(t)$; From
\eqref{w.corretto}, \eqref{def.a2} we see that $R_{1}^{1 }(t)$ can be written
as follows:
\begin{equation}\label{formula.R2}
  \begin{split}
  &  R_{1}^{1 }(t):=   J_1(t)+ K_{1}^{1 }(t) \, ,  \quad \textrm{where } \,  J_1(t)=   \eu^{Lt} (\Id- \PP_1) x(0) \, ,  \quad  \textrm{and} \\
  & K_{1}^{1 }(t):=   \int_0^t \eu^{L (t-s)} (\Id-\PP_1) M_{1}^{1,R }(s)ds -\int_t^{+\infty} \eu^{L (t-s)}\PP_1 M_{1}^{1,R }(s) ds   \\
  &  M_{1}^{1,R }(s)=   N(x(s))-N(A_{1}^{0 }(s))= N(A_{1}^{1 }(s)+R_{1}^{1 }(s))-N(A_{1}^{0 }(s)).
  \end{split}
\end{equation}
 We stress that  by construction $J_1(t)= o(\eu^{(-\wedge_2+\ep)s})$.
Since $R_{1}^{1 }(t)=R_{1}^{0 }(t)-a_{1}^{1 }(t)=o(\eu^{-(2\wedge_1-\ep)t })$,
cf. \eqref{def.R1},\eqref{asymp.a2}, we find that $M_{2,R}(t)=O(A_{1}^{0 }(t)R_{1}^{1 }(t))=o(
\eu^{-(3\wedge_1-\ep)t})$.  Therefore, from \eqref{formula.R2} and  Lemma \ref{stimazero},  we find
$K_{1}^{1 }(t) =o(\eu^{-(3\wedge_1-\ep)t })$ so that $R_{1}^{1 }(t) =o(\eu^{
  -(\wedge_2-\ep)t}+\eu^{-(3\wedge_1-\ep)t })$.
\smallskip

So we can expand $x(t)$ as
$x(t)= A_{1}^{1 }(t)+R_{1}^{1 }(t)$ and $R_{1}^{1 }(t)=o(\eu^{
  -(\wedge_2-\ep)t}+\eu^{-(3\wedge_1-\ep)t })$.
If  $k_2=2$, then  $R_{1}^{1 }(t) =o(\eu^{ -(\wedge_2-\ep)t})$ and
Lemma \ref{l:step.a} is proved, while if
 $k_2>1$, then $R_{1}^{1 }(t) =o(\eu^{-(3\wedge_1-\ep)t })$ and  we go
 to the next step.
\smallskip

\noindent \textbf{--Step 2. } Since
$N(A_{1}^{1 }(t))-N(A_{1}^{0 }(t))=o( \eu^{-(3\wedge_1-\ep)t })$, we
can define
\begin{equation}\label{def.a3}
  \begin{split}
&  M_{1}^{2 }(s):=N(A_{1}^{1 }(s))-N(A_{1}^{0 }(s))    \, , \quad  \textrm{and} \\
&  a_{1}^{2 }(t):=  \int_{0}^{t}  \eu^{L(t-s)}(\Id-\PP_1) M_{1}^{2 }(s) ds - \int_{t}^{+\infty} \eu^{L(t-s)}  \PP_1 M_{1}^{2 }(s) ds.
  \end{split}
\end{equation}
Notice that
\begin{equation}\label{asymp.a3}
  a_{1}^{2 }(t)=o(\eu^{-(3\wedge_1-\ep)t })\,\, \textrm{ as $t \to +\infty$. }
\end{equation}
Denote by $A_{1}^{2 }(t)=A_{1}^{1 }(t)+a_{1}^{2 }(t)$, and by $R_{1}^{2 }(t)=x(t)-A_{1}^{2 }(t)$. From
\eqref{w.corretto}, \eqref{def.a2}, \eqref{def.a3} we see that $R_{1}^{2 }(t)$
solves the following fixed point problem:
\begin{equation}\label{formula.R3}
  \begin{split}
    R_{1}^{2 }(t):= & J_1(t)+ K_{1}^{2 }(t) \qquad \qquad  \textrm{where}
   \\
   K_{1}^{2 }(t)= &  \int_0^t \eu^{L (t-s)} (\Id-\PP_1) M_{1}^{2,R }(s)ds-\int_t^{+\infty} \eu^{L (t-s)}\PP_1 M_{1}^{2,R }( s) ds   \\
    M_{1}^{2,R }( s)=& N(x(s))- N(A_{1}^{1 }(s))=N(A_{1}^{2 }(s)+R_{1}^{2 }(s))-N(A_{1}^{1 }(s))
  \end{split}
\end{equation}

Since $R_{1}^{2 }(t)=R_{1}^{1 }(t)-a_{1}^{2 }(t)=o(\eu^{-(3\wedge_1-\ep)t })$,
cf. \eqref{asymp.a3}, we find that $M_{1}^{2,R }(t)=O(A_{1}^{0 }(t)R_{1}^{2 }(t))=o(
\eu^{-(4\wedge_1-\ep)t})$. Therefore from \eqref{formula.R3} and  Lemma \ref{stimazero} we
find that  $K_{1}^{2 }(t) =o( \eu^{-(4\wedge_1-\ep)t })$.

So we can expand $x(t)$ as
$x(t)= A_{1}^{2 }(t)+R_{1}^{2 }(t)$ and $R_{1}^{2 }(t)=o(\eu^{ -(\wedge_2-\ep)t}+ \eu^{-(4\wedge_1-\ep)t })$. If  $k_2=2$ then
$R_{1}^{2 }(t) =o(\eu^{ -(2\wedge_2-\ep)t})$ and Lemma \ref{l:step.a} is proved, while if
 $k_2>2$ then $R_{1}^{2 }(t) =o(\eu^{-(4\wedge_1-\ep)t })$ we
 iterate the argument of Step j below, till $j=k_2$.
 \smallskip

\noindent \textbf{--Step j. } Since
$N(A_1^{j-1 }(t))-N(A_1^{j-2 }(t))=o(\eu^{ -[(j+1)\wedge_1-\ep]t })$ as $t \to
+\infty$, we can define
\begin{equation}\label{def.aj}
  \begin{split}
&  M_1^{j }(s):=N(A_1^{j-1 }(s))-N(A_1^{j-2 }(s))    \, , \qquad \textrm{and} \\
&  a_{1}^{j }(t):=  \int_{0}^{t}  \eu^{L(t-s)}(\Id-\PP_1) M_{1}^{j }(s) ds  - \int_{t}^{+\infty} \eu^{L(t-s)}  \PP_1 M_{1}^{j }(s) ds
  \end{split}
\end{equation}
Notice that
\begin{equation}\label{asymp.aj}
a_{1}^{j }(t)=o(\eu^{-[(j+1)\wedge_1-\ep]t })\,\, \textrm{ as $t \to +\infty$. }
\end{equation}
Denote by $A_{1}^{j }(t)=A_1^{j-1 }(t)+a_1^{j }(t)$, and by
$R_1^{j }(t)=x(t)-A_{1}^{j }(t)$. From \eqref{w.corretto}, \eqref{def.a2},
\eqref{def.a3}, \eqref{def.aj} we see that $R_1^{j }(t)$ can be written
as follows:
\begin{equation}\label{formula.Rj}
  \begin{aligned}
   & R_1^{j }(t):=  J_1(t)+ K_1^{j }(t) \quad \textrm{where}
   \\
   &K_1^{j }(t)=   \int_0^t \eu^{L (t-s)} (\Id-\PP_1) M_1^{j,R }(s)ds-\int_t^{+\infty} \eu^{L (t-s)}\PP_1 M_1^{j,R }( s) ds   \\
   & M_1^{j,R }( s)= N(x(s))-N(A_{1}^{j-1 }(s))=N(A_{1}^{j }(s)+R_1^{j }(s))-N(A_{1}^{j-1 }(s)).
  \end{aligned}
\end{equation}
Since $R_1^{j }(t)=R_1^{j-1 }(t)-a_1^{j }(t)= o(\eu^{-[(j+1)\wedge_1-\ep]t })$,
cf \eqref{asymp.aj}, we find that the term
$M_1^{j,R }( t)=O(A_{1}^{0 }(t)R_1^{j }(t))=o(
\eu^{-[(j+2)\wedge_1-\ep]t})$. Thereafter, from \eqref{formula.Rj} and Lemma \ref{stimazero}
we find     $K_1^{j }(t)=o(\eu^{-[(j+1)\wedge_1-\ep]t })$, and
 $R_1^{j }(t) =o(\eu^{-[(j+1)\wedge_1-\ep]t }+ \eu^{-(\wedge_2-\ep) t})$.
\end{proof}
In the next table we enumerate the terms $a_i^j(t)$, and the related
asymptotic behaviors,  that we are going to use in
the remaining part of the proof.
\begin{center}
\begin{longtable}{|c|c|c|cc|cc|}
\hline
$ \underset{e^{-\wedge_1t}}{\quad {\it{a}_1^0}(t)}$ &
$\underset{e^{-2\wedge_1t}}{\quad {\it{a}_1^1}(t)}$ &
$\underset{e^{-3\wedge_1t}}{\quad {\it{a}_1^2}(t)}$ &
$\ldots$ &  & $\underset{e^{-(k_2+1)\wedge_1t}}{\quad {\it{a}_1^{k_2}}(t)} $ &\\
\hline
$ \underset{e^{-\wedge_2t}}{\quad {\it{a}_2^0}(t)} $ &
$\underset{e^{-(\wedge_1 +\wedge_2)t}}{\quad {\it{a}_2^1}(t)} $  &
$\underset{e^{-(2\wedge_1 +\wedge_2)t}}{\quad {\it{a}_2^2}(t)} $ &
$\underset{e^{-(3\wedge_1 +\wedge_2)t}}{\quad {\it{a}_2^3}(t)} $&
$\ldots$ &
$\underset{e^{-(k_3\wedge_1 +\wedge_2)t}}{\quad {\it{a}_2^{k_3}}(t)} $ &\\
\hline
$\underset{e^{-\wedge_3 t}}{\quad {\it{a}_3^0}(t)} $ &
$\underset{e^{-(\wedge_1 +\wedge_3)t}}{\quad {\it{a}_3^1}(t)}$ &
$\underset{e^{-(2\wedge_1 +\wedge_3)t}}{\quad {\it{a}_3^2}(t)} $ &
$\underset{e^{-(3\wedge_1 +\wedge_3)t}}{\quad {\it{a}_3^3}(t)} $ &
$\ldots$ &
$\underset{e^{-(k_4\wedge_1
    +\wedge_3)t}}{\quad {\it{a}_3^{k_4}}(t)} $ &\\
\hline
$\vdots$ & $\vdots$& $\vdots$& $\vdots$& $\vdots$&$\vdots$ &\\
\hline
$\underset{e^{-\wedge_m t}}{\quad {\it{a}_m^0}(t)} $ &
$\underset{e^{-(\wedge_1 +\wedge_m)t}}{\quad {\it{a}_m^1}(t)} $ &
$\underset{e^{-(2\wedge_1 +\wedge_m)t}}{\quad {\it{a}_m^2}(t)} $ &
$\underset{e^{-(3\wedge_1 +\wedge_m)t}}{\quad {\it{a}_m^3}(t)} $ &
$\ldots$ &
$\underset{e^{-(r\wedge_1
    +\wedge_m)t}}{\quad {\it{a}_m^{r}}(t)} $ & $\ldots$
    \\
\hline
    \caption{ 
   \footnotesize Terms $a_i^j(t)$:  We omit the $\ep$-shift of the exponents so that
$\eu^{-\wedge t}$ stands for $o(\eu^{-(\wedge-\ep) t})$.} \label{table-1}
\end{longtable}
\end{center} \vspace{-1 cm}
In the next table   we list  the terms $R_i^j(t)$, and the related
asymptotic behaviors.
\begin{center}
\begin{longtable}{|c|c|c|cc|cc|}
\hline
$ \underset{e^{-2\wedge_1t}}{\quad {\it{R}_1^0}(t)}$ &
$\underset{e^{-3\wedge_1t}}{\quad {\it{R}_1^1}(t)}$ &
$\ldots$ &
$\underset{e^{-(k_2+k_1)\wedge_1t}}{\quad {\it{R}_1^{k_2-1}}(t)}$ &  & $\underset {e^{-\wedge_2t}}{\quad {\it{R}_1^{k_2}}(t)}$ &\\
\hline
$ \underset{e^{-(k_2+k_1+1) \wedge_1 t}}{\quad {\it{R}_2^0}(t)}$ &
$\underset{e^{-(k_2+k_1+2)\wedge_1  t}}{\quad {\it{R}_2^1}(t)} $  &
$\ldots$ &
$\underset{e^{-(k_3+k_2+k_1)\wedge_1t}}{\quad {\it{R}_2^{k_3-1}}(t)}  $&
  &
$\underset{e^{-\wedge_3 t}}{\quad {\it{R}_2^{k_3}}(t)} $ &\\
\hline
$\underset{e^{-(k_3+k_2+k_1+1)\wedge_1 t}}{\quad {\it{R}_3^0}(t)} $ &
$\underset{e^{-(k_3+k_2+k_1+2)\wedge_1t}}{\quad {\it{R}_3^1}(t)} $ &
$\ldots$ &
$\underset{e^{-(k_4+k_3+k_2+k_1)\wedge_1 t}}{\quad {\it{R}_3^{k_4-1}}(t)} $ &
  &
$\underset{e^{-\wedge_4 t}}{\quad {\it{R}_3^{k_4}}(t)} $ &\\
\hline
$\vdots$ & $\vdots$& $\vdots$& $\vdots$&  &$\vdots$ &\\
\hline
$\underset{e^{-(k_{m}+\ldots+k_{1}+1) \wedge_1 t}}{\quad {\it{R}_m^0}(t)} $ &
$\underset{e^{-(k_{m}+\ldots+k_{1}+2) \wedge_1 t}}{\quad {\it{R}_m^1}(t)} $ &
$\ldots$ &
$\underset{e^{-(k_{m}+\ldots+k_{1}+r+1)\wedge_1 t}}{\quad {\it{R}_m^{\nu}}(t)}  $ &
 &  $\ldots$ &
  \\
\hline
    \caption{
\footnotesize Terms $R_i^j(t)$:  Also in this case we omit the $\ep$-shift of the exponents so that
$\eu^{-\wedge t}$ stands for $o(\eu^{-(\wedge-\ep) t})$.}  \label{table-2}
\end{longtable}
\end{center} \vspace{-1 cm}

It is worthwhile to observe that $a_{1}^{1 }(t)$ and  $a^{j}_1(t)$ are
solutions respectively of the (explicitly solvable) linear
non-homogeneous  problems
\begin{equation}\label{step2}
    \dot{y}= Ly +M_{1}^{1 }(t), \,\,\textrm{ and } \,\,
    \dot{y}= Ly +M_{1}^{j }(t),
\end{equation}
where $M_{1}^{1 }(t)$ and $M_1^j(t)$ are defined   in
\eqref{def.a2} and in \eqref{def.aj},
respectively. The initial conditions can be written respectively as
\[
a_{1}^{1 }(0)=-\int_0^{+\infty} \eu^{-Ls}\PP_1 M_{1}^{1 }(s) ds,\,\,
 \textrm{ and } \,\, a_1^j(0)=-\int_0^{+\infty} \eu^{-Ls}\PP_1
 M_1^j(s) ds.
\]
 However $a_{1}^{1 }(0)$, $a_1^j(0)$   are the real unknown of the problem,
since they depend on $\bar{\ell}_1$, which is evaluated by using a
fixed point argument and not by a closed formula. In fact  $\bar{\ell}_1$  can just
be approximated as $\bar{\ell}_1=\PP_1 x(0)+O(x^2(0))$.

  Also, the remainder
terms $K_{1}^{1 }(t)$ and $K_1^{j }(t)$  solve the non-linear differential
equations
\begin{equation}\label{eq.R2}
    \dot{y}= Ly +M_1^{R,1}( t),\,\, \textrm{ and } \,\,
    \dot{y}= Ly +M_1^{R,j}( t),
\end{equation}
where $M_1^{R,1}( t)$ and $M_1^{R,j}( t)$ are defined   in
\eqref{formula.R2} and in \eqref{formula.Rj}, respectively (and are non linear function of $y$).

In order to proceed with the expansion beyond $\eu^{-\wedge_2 t}$ we need to take into account the contribution of the linear part once again.
For this reason we set
  \begin{equation}\label{N2}
\bar{\ell}_2= P_2 (x(0))+ \bar{\mathfrak{N}}_2 ,
 \,\, \textrm{ where } \,\,
   \bar{\mathfrak{N}}_2:= \int_{0}^{+\infty} \eu^{-Ls} P_2 [N(x(s))- N(A_1^{k_2}(s))] ds.
\end{equation}
Notice that $ \bar{\mathfrak{N}}_2$ is well defined since
\begin{equation}\label{NN}
M^{0,R}_2(t):=N(x(t))- N(A_1^{k_2}(t))=o(\eu^{-[(k_{2}+k_1+1)\wedge_1-\ep]t}) \,.
\end{equation}
 Then, we iterate the previous argument, and we prove the following.
\begin{lemma}\label{l:step.b}
  Let us recall that   $ \wedge_1 \left(\sum_{i=1}^{3}k_i \right) \le \wedge_3 <   \wedge_1\left(1+\sum_{i=1}^{3}k_i \right)$, cf \eqref{defkappa}, and set
$a_2^{0 }(t):=\ell_2(t)= \eu^{L t}  \bar{\ell}_2$. Then we can expand $x(t)$ as follows
 \begin{equation}\label{expand.b.fin}
 \begin{split}
 &  x(t)=A_1^{k_{2}}(t)+ \sum_{j=0}^{ k_{3}}  a^j_2(t) +R_2^{k_{3}}(t) \, , \quad \textrm{where}  \\
 & a^j_2(t)=o(\eu^{-(\wedge_2+j \wedge_1-\ep)t}) \, , \quad \textrm{and }  \; R_2^{k_3}(t)=o(\eu^{-  (\wedge_3-\ep)t})
 \end{split}
 \end{equation}
\end{lemma}
\begin{proof}
Let us set
$A_2^{-1}(t)= A_1^{k_2}(t)$, $A_2^{0 }(t)=A_2^{-1}(t)+a_2^{0 }(t)$.
From \eqref{w.nonlin}, \eqref{formula.Rj}, we get
\begin{equation*}
  \begin{split}
        R_2^{0}(t)& =x(t)-A_2^{0 }(t)= \eu^{L t} (\Id- \mathcal{P}_1) x(0)- \eu^{L t}  P_2 x(0) \\ &+
        \int_{0}^{t} \eu^{L (t-s)} (\Id -P_1-P_2) M^{0,R}_2(s) ds
        - \int_{t}^{+\infty} \eu^{L (t-s)} (P_1+P_2) M^{0,R}_2(s) ds
  \end{split}
\end{equation*}
and consequently
\begin{equation}\label{J2}
  \begin{split}
       R_2^{0}(t)= &x(t)-A_2^{0 }(t)= J_2(t)+K^0_2(t) \, , \quad
       \textrm{ where } \,
 J_2(t)=\eu^{L t} (\Id- \mathcal{P}_2) x(0)  \\
         K_2^0(t)= &\int_{0}^{t} \eu^{L (t-s)} (\Id -\mathcal{P}_2)
         M^{0,R}_2(s) ds -
\int_{t}^{+\infty} \eu^{L (t-s)} \mathcal{P}_2 M^{0,R}_2(s) ds.
  \end{split}
\end{equation}
Notice that $a_2^{0 }(t)=o(\eu^{-(\wedge_2-\ep)t})$, $J_2(t)=o(\eu^{-(\wedge_3-\ep)t})$. Further,
using Lemma \ref{stimazero} and \eqref{NN}, we find $K_2^0(t)=o(\eu^{-[(k_2+k_1+1)\wedge_1-\ep]t})$
and
\begin{equation}\label{est.R0B}
R_2^{0}(t)=o(\eu^{-[(k_2+k_1+1)\wedge_1-\ep]t}+\eu^{-(\wedge_3-\ep)t}) \, .
\end{equation}
Whence, if $k_3=0$, then $R_2^{0}(t)=o(\eu^{-(\wedge_3-\ep)t})$, and the Lemma is proved.

If $k_3>0$, we proceed  as above; We denote by
\begin{equation}\label{def.b1}
  \begin{split}
&  M_2^1(s):=N(A_2^{0 }(s))-N(A_2^{-1}(s))    \, , \qquad \textrm{and} \\
&  a^{1}_2(t):=  \int_{0}^{t}  \eu^{L(t-s)}(\Id-\PP_2) M_2^1(s) ds  - \int_{t}^{+\infty} \eu^{L(t-s)}  \PP_2  M_2^1(s) ds
  \end{split}
\end{equation}
and in general for any $1 \le j \le k_3$ we set
\begin{equation}\label{def.bj}
  \begin{split}
&  M_2^{j}(s):=N(A_2^{j-1}(s))-N(A_2^{j-2}(s))    \, , \qquad \textrm{and} \\
&  a_2^{j}(t):=  \int_{0}^{t}  \eu^{L(t-s)}(\Id-\PP_2) M_2^{j}(s) ds
-
 \int_{t}^{+\infty} \eu^{L(t-s)}  \PP_2  M_2^{j}(s) ds.
  \end{split}
\end{equation}
Then, using again Lemma \ref{stimazero},  we  find that
\begin{equation}\label{est.b1}
  \begin{split}
&  M_2^1(t) =O\big( a_2^{0}(t)   \cdot  A_2^{-1}(t)\big) =O \big(  a_2^{0}(t)  \cdot   a_{1}^{0 }(t) \big)   =o\big( \eu^{-(\wedge_2 +\wedge_1-\ep)t}\big)\\
&  a_2^{1}(t) =o\big( \eu^{-(\wedge_2 +\wedge_1-\ep)t}\big)
  \end{split}
\end{equation}
Analogously, through an inductive argument, for any $j \in \{1, \ldots, k_3 \}$ we get
\begin{equation*}
  \begin{split}
&  M_2^{j}(t) =O\big( a_2^{j-1}(t)   \cdot  A_2^{j-2}(t)\big) =O \big(  a_2^{j-1}(t)  \cdot   a_{1}^{0 }(t) \big)   =o\big( \eu^{-[\wedge_2 +j\wedge_1-\ep]t}\big)\\
&  a_2^{j}(t) =o\big( \eu^{-[\wedge_2 +j\wedge_1-\ep]t}\big).
  \end{split}
\end{equation*}
As above  the remainder term  $R_{2}^{1}(t):=x(t)-A_2^1(t)$ can be
written as follows
\begin{equation}\label{formula.R1B}
  \begin{split}
    R^{1}_2(t):= &   J_2(t)+ K^1_2(t) \, ,  \quad \textrm{where } \\
    K^1_2(t):= &  \int_0^t \eu^{L (t-s)} (\Id-\PP_2) M^{1,R}_2(s)ds -\int_t^{+\infty} \eu^{L (t-s)}\PP_2 M^{1,R}_2(s) ds   \\
     M^{1,R}_2( s)= & N(x(s)-N(A_2^{0 }(s))= N(A_2^1(s)+R^{1}_2(s))-N(A_2^{0 }(s))
  \end{split}
\end{equation}
Since $R^{1}_2(t)=R^{0}_2(t)-a_2^1(t)  =o(\eu^{-[(k_{2}+k_1+1)\wedge_1-\ep]t })$, see \eqref{est.R0B}, \eqref{est.b1}, we find
 \begin{equation*}
M^{1,R}_2(
t)=O(A_2^{0 }(t)R^{1}_2(t))=O(a_{1}^{0 }(t)R^{1}_2(t))=o(\eu^{-[(k_{2}+k_1+2)\wedge_1-\ep]t
}).
\end{equation*}
Hence, from Lemma \ref{stimazero}
 we find $K^1_2(t)=o(\eu^{-[(k_2+k_{1}+2)\wedge_1-\ep]t })$.
 So we get
 \begin{equation}\label{est.R1B}
 R^{1}_2(t)=o(\eu^{-[(k_2+k_{1}+2)\wedge_1-\ep]t }+\eu^{-(\wedge_3-\ep)t }) \,.
\end{equation}
Whence, if $k_3=1$, then $R_2^1(t)=o(\eu^{-(\wedge_3-\ep)t})$, and the Lemma is proved.

Otherwise we proceed by induction assuming that $R^{j-1}_2(t)
=o(\eu^{-[(k_2+k_{1}+j)\wedge_1-\ep]t }+\eu^{-(\wedge_3-\ep)t })$
 and proving that
 $R^{j}_2(t)=o(\eu^{-[(k_2+k_{1}+j+1)\wedge_1-\ep]t }
+\eu^{-(\wedge_3-\ep)t })$, for $1 \le j \le k_3$. Note that
 \begin{equation}\label{formula.R2B}
  \begin{split}
    R^{j}_2(t):= &   J_2(t)+ K^{j}_2(t) \, ,  \quad \textrm{where } \\
    K^{j}_2(t):= &  \int_0^t \eu^{L (t-s)} (\Id-\PP_2) M^{j,R}_2(s)ds -\int_t^{+\infty} \eu^{L (t-s)}\PP_2 M^{j,R}_2(s) ds   \\
     M^{j,R}_2(s)= & N(x(s)-N(A_2^{j-1}(s))= N(A_2^{j}(s)+R^{j}_2(s))-N(A_2^{j-1}(s))
  \end{split}
\end{equation}
 In fact $R^{j}_2(t)=R^{j-1}_2(t)-a_2^{j}(t)=o(\eu^{-[(k_2+k_{1}+j)\wedge_1-\ep]t})$, whence
 \begin{equation*}
   M^{j,R}_2(t)=
   O[R^{j}_2(t)A_2^j(t)]=O[R^{j}_2(t)a_{1}^{0 }(t)]=o(\eu^{-[(k_2+k_{1}+j+1)\wedge_1-\ep]t
   }).
\end{equation*}
 Therefore,  from Lemma \ref{stimazero} we find
 $R^{j+1}_2(t)=o(\eu^{-[(k_2+k_{1}+j+1)\wedge_1-\ep]t }+\eu^{- (\wedge_3-\ep)t })$, so the claim is proved. In particular,
  iterating the argument till $j+1=k_3$
  we find
  $R^{\kappa_3}_2(t)=o(\eu^{- (\wedge_3-\ep)t })$, and we conclude the proof
 of the Lemma.
 \end{proof}
\begin{remark}\label{Reb}
  We stress that $R_2^{j}(t)$ is in general negligible with respect to $a_2^j(t)$ for any $j$.
   In fact
  $a_2^j(t)=o(\eu^{-(\wedge_2 +j \wedge_1-\ep)t})$,  $R_2^{j}(t)=o(\eu^{-[(k_2+k_{1}+1) \wedge_1-\ep]t})$, and
  $\wedge_2 +j \wedge_1 < (k_2+k_{1}+1)\wedge_1+ j\wedge_1$ .
\end{remark}
Notice that we can continue this iterative scheme up to an arbitrary order. In fact,
 recalling the definitions of $a^j_1(t)$, $A^j_1(t)$, $a_2^j(t)$, $A^j_2(t)$ of Lemmas \ref{l:step.a}, \ref{l:step.b}, and
 the definitions of the $k_i$ given in \eqref{defkappa}, we set
$$  A_{i}^{\mu}(t)=A_{i-1}^{k_{i}}(t)+ \sum_{j=0}^{\mu} a_i^{j}(t) \, , \quad \textrm{for any $0 \le \mu \le k_{i+1}, \;$ $i=2, \ldots, m$} $$
  (setting also $k_{m+1}=+\infty$), where the functions $a_i^j(t)$ are defined as follows
\begin{equation}\label{redefine1}\begin{split}
a_i^{0}(t):= & \, \eu^{L t} [ P_i x(0)+ \bar{\mathfrak{N}}_i] \, , \quad \textrm{for $i=2, \ldots , m \quad$ where}\\
  \bar{\mathfrak{N}}_i:= & \int_{0}^{+\infty} \eu^{-Ls} P_i [N(x(s))- N(A_{i-1}^{k_i-1}(s))] ds \, .
 \end{split}
\end{equation}
If $k_i>0$ we also have
\begin{equation}\label{redefine2}\begin{split}
 M_i^{1}(s):= & \, N(A^{0}_i(s))-N(A^{k_i}_{i-1}(s))  \, ; \\
a_i^{1}(t):=  & \int_{0}^{t}  \eu^{L(t-s)}(\Id-\PP_i) M_i^{1}(s) ds
- \int_{t}^{+\infty} \eu^{L(t-s)}  \PP_i M_i^{1}(s) ds \, .
\end{split}
\end{equation}
If $k_i \ge j>1$ we also have
\begin{equation}\label{redefine3}\begin{split}
M_i^{j}(s):= & N(A^{j-1}_i(s))-N(A^{j-2}_i(s))   \quad \textrm{for $1 \le j \le k_{i+1}$   }\\
a_i^{j}(t):=  & \int_{0}^{t}  \eu^{L(t-s)}(\Id-\PP_i) M_i^{j}(s) ds  - \int_{t}^{+\infty} \eu^{L(t-s)}  \PP_i  M_i^{j}(s) ds.
\end{split}
\end{equation}

We have the following result,  see also Table~\ref{table-1} and Table~\ref{table-2}.

\begin{lemma}\label{expand.finale}
  Assume $N \in C^2$ and let $x(t)$ be a solution of \eqref{sistemacompleto} such that $x(t) \to 0$ as $t \to +\infty$.
    Let $r \in \mathbb{N}$, and set $\mu=r + \sum_{i=1}^m k_i$, then we can expand $x(t)$ as follows
 \begin{equation}\label{expand.b.fin}
 \begin{split}
 &  x(t)=\sum_{j=0}^{k_{2}} a^{j}_1(t)+ \ldots + \sum_{j=0}^{k_m}a^{j}_{i-1}(t)+ \sum_{j=0}^{\mu} a^j_m(t)+R^{\mu}_m(t)
  \, , \quad \textrm{where}  \\
 & a^j_i(t)=o(\eu^{-(\wedge_i+j \wedge_1-\ep)t}) \, , \quad \textrm{and }  \; R^{\mu}_m(t)=o(\eu^{-  (\mu+1-\ep) \wedge_1 t})
 \end{split}
 \end{equation}
\end{lemma}
\noindent The proof is simply an iteration of the previous scheme.

\medskip

Now, using the fact that all the functions appearing in the definition of $a^j_i(t)$ are exponentials, possibly  multiplied by polynomials,
 we can improve the estimates of Lemma \ref{expand.finale}.
To proceed we need more regularity in order to replace the function $N$ of \eqref{sistemacompleto} by its Taylor polynomial of degree $\mu \ge 2$, say $N_{\mu}$.
Let us fix $r \in \mathbb{N}$ and $\mu =r+ \sum_{i=1}^m k_i$, and go back to Step 1.
We rename $\tilde{a}^0_1(t):=a^0_1(t)$, see \eqref{def.R1}, then we
  denote by $\tilde{a}^1_1(t)$  the function
$a_{1}^{1 }(t)$ of \eqref{def.a2} but where $N$ is replaced by $N_{\mu}$. Notice that $\tilde{a}^1_1(t)-a_{1}^{1 }(t)=o(\eu^{-(\mu+1-\ep) \wedge_1 t})$
so that $\tilde{R}^1_1(t)=x(t)-\tilde{a}^1_1(t) \sim R_{1}^{1 }(t)$. Then we iterate the scheme above replacing everywhere $N$ by $N_{\mu}$
and obtaining new functions denoted by `` $\tilde{\cdot}$ ".
I.e. we set $\tilde{a}^0_1(t)=a^0_1(t)$,
\begin{equation*}
\begin{split}
     \tilde{A}_1^{r}(t)=  & \sum_{j=0}^{r} \tilde{a}_1^j(t) \; \textrm{ for $2 \le r \le k_{2}$} \, , \qquad
  \tilde{A}_2^{r}(t)= \tilde{A}_1^{k_2}(t)+ \sum_{j=0}^{r}\tilde{a}_2^j(t)\; \textrm{ for $0 \le r \le k_3$}\, , \\
     \tilde{A}_i^{r}(t)=  & \tilde{A}_{i-1}^{k_{i}}(t)+ \sum_{j=0}^{r} \tilde{a}_i^j(t)\; \textrm{ for $0 \le r   \le k_{i+1}$}
  \end{split}
\end{equation*}
  where $\tilde{\mathfrak{N}}_i$,   $\tilde{M}^j_i(t)$, $\tilde{a}^j_i(t)$ are defined as
  $\bar{\mathfrak{N}}_i$,   $M^j_i(t)$, $a^j_i(t)$ in \eqref{redefine1}, \eqref{redefine2}, \eqref{redefine3} but with $N(\cdot)$, $A^j_i(t)$ replaced by $N_{\mu}(\cdot)$ and $\tilde{A}^j_i(t)$.
  Then we have the following.
\begin{lemma}\label{expand.finalebis}
 Let $r \in \mathbb{N}$,  set $\mu=r + \sum_{i=1}^m k_i$, and assume $N \in C^{\mu+1}$;
 denote by $N_{\mu}(x)$ the Taylor polynomial of degree $\mu$ centered in the origin.
  Let $x(t)$ be a solution of \eqref{sistemacompleto} such that $x(t) \to 0$ as $t \to +\infty$, then
  set $\tilde{\mathfrak{N}}_1:= \bar{\mathfrak{N}}_1$, see \eqref{N0}, and
  \begin{equation}\label{redefinebis}\begin{split}
 &  \tilde{\mathfrak{N}}_i:=  \int_{0}^{+\infty} \eu^{-Ls} P_i [N (x(s))- N_{\mu}(\tilde{A}_{i-1}^{\kappa_{i}-1}(s))] ds
   \quad \textrm{for $i=2, \ldots , m$}\\
& \tilde{a}_i^{0}(t):=  \eu^{L t} [ P_i x(0)+ \tilde{\mathfrak{N}}_i] \, , \quad \textrm{for $i=1, \ldots , m $.}
\end{split}
\end{equation}
Then we set $\tilde{A}^{-1}_i(s):= \tilde{A}^{k_i}_{i-1}(s)$ for $i=2, \ldots, m$, and for $1 \le j \le k_{i+1}$
  \begin{equation}\label{redefinetris}
  \begin{split}
\tilde{M}_i^{j}(s):= & N_{\mu}(\tilde{A}^{j-1}_i(s))-N_{\mu}(\tilde{A}^{j-2}_i(s))\\
\tilde{a}_i^{j}(t):=  & \int_{0}^{t}  \eu^{L(t-s)}(\Id-\PP_i) \tilde{M}_i^{j}(s) ds  - \int_{t}^{+\infty} \eu^{L(t-s)}  \PP_i  \tilde{M}_i^{j}(s) ds
\end{split}
\end{equation}
 ($k_{m+1}=+\infty$). Then we can expand $x(t)$ as follows:
 \begin{equation}\label{expand.b.finbis}
 \begin{split}
 &  x(t)=\tilde{A}^{\mu}_m(t)+ \tilde{R}^{r}_m(t)
 \end{split}
 \end{equation}
 where again $\tilde{R}_m^{r}(t)=o(\eu^{-  [(\mu+1)\wedge_1-\ep]t})= o(\eu^{-  [(1+r+\sum_{i=1}^m k_i)\wedge_1-\ep]t}) $ but $\tilde{A}^{\mu}_m(t)$ is as above.
\end{lemma}
\begin{proof}
Formula \eqref{expand.b.finbis} simply follows by observing that, at each step, $\tilde{A}^j_i(t)-A^j_i(t)=o(\eu^{-(\mu+1-\ep) \wedge_1 t})$
so that
 $\tilde{R}_i^j(t)=x(t)- \tilde{A}_i^j(t)$ and $R_i^j(t)=x(t)- A_i^j(t)$  have the same order.
\end{proof}
Our purpose now is to observe that each coordinate of $\tilde{A}^{\mu}_m(t)$ is made up by sum of exponentials, possibly multiplied by polynomials, in the resonant
cases.
Therefore we introduce the following assumptions:
\begin{description}
  \item[$\bs{R1}$] All the eigenvalues $\la_i$ of $L$ are real and distinct, hence  $k_i=1$ for any $i$, and $m=n$.
  \item[$\bs{R2}$] $\bs{R1}$ holds  and for any $i \in 1, \ldots , n$ and $(\chi_1, \ldots, \chi_{i-1}) \in \mathbb{N}^{i-1}$,
  then  $\wedge_i \ne \sum_{j=1}^{i-1} \wedge_j \chi_j$.
\end{description}
We assume first $\bs{R2}$ so that
  our iterative scheme contains no resonances at all.

 \begin{lemma}\label{nonres}
  Let $r \in \mathbb{N}$, and $\mu=r + \sum_{i=1}^m k_i$; Assume $N \in C^{\mu+1}$, and  $\bs{R2}$; then
  the function $\tilde{A}^{\mu}_m(t)$ of Lemma \ref{expand.finalebis} is a vectorial polynomial of degree (at most) $\mu$, i.e.
   each coordinate  $\tilde{A}^{\mu}_{m,s}(t)$ of
  $\tilde{A}^{\mu}_m(t)$   is a polynomial
 in the functions $\eu^{-\wedge_i t}$. More precisely,
 there is a polynomial $p_s^{\mu}(f_1, \ldots , f_n)$ of degree $\mu$ such that $$\tilde{A}^{\mu}_{m,s}(t)=p_s^{\mu}(\eu^{-\wedge_1 t}, \ldots , \eu^{-\wedge_n t}) \, .$$
 \end{lemma}
 \begin{proof}
  Requiring $\bs{R1}$ we see that there is $c_{i,s} \in \RR$ such that the $s^{th}$ coordinate of the linear terms $\tilde{a}^{0}_i(t)$
is $c_{i,s} \eu^{-\wedge_i t}$ while the others are null
(we do not have linear resonances).
 $\bs{R2}$ is needed to   avoid also nonlinear resonances, i.e. resonances in the integral \eqref{redefinetris} defining $\tilde{a}^j_i(t)$
 for $j \ge 1$.
 In fact it is easy to check that the first coordinate of   $\tilde{M}^1_1(s)$ defined in \eqref{redefinetris} (or in \eqref{def.a2}  but with $N_{\mu}$ replacing $N$) is a polynomial
 in $\eu^{-\wedge_1 t}$. Then   $\bs{R2}$ guarantees that each coordinate $s$ of $\tilde{a}^1_1(t)$ takes the form
 $p_s^{1}(\eu^{-\wedge_1 t}, \ldots , \eu^{-\wedge_n t})$  where $p_s^{1}$ is a polynomial of degree $\mu$. Then we go through all the iterative scheme
 and we get the thesis of the Lemma.
 \end{proof}
 Let us introduce some notation in  line with Section 2.4.
  \begin{equation}\label{Imu}
    I_\Theta=\big\{\chi=(\chi_1, \ldots ,\chi_m) \in \mathbb{N}^m \, \colon \,  \chi_1
    \wedge_1+ \ldots  +\chi_m \wedge_m \le \Theta \big\}.
  \end{equation}
  and observe that if $\bs{R1}$ holds then  $m=n$.

 Now we reformulate the previous Lemma in a way that is more suitable to be used for the parabolic problem studied in this article.
   \begin{prop}\label{nonres2}
 Assume $N \in C^{\infty}$,   $\bs{R2}$ and let $x(t)=(x_1(t), \ldots, x_n(t)) \to 0$ as $t \to +\infty$; then for any $\Theta>0$ we can find   coefficients
 $c^\chi=(c^\chi_{1}, \ldots, c^\chi_{n})\in \RR^n$, $\chi \in I_{\Theta}$ such that for any $s=1, \ldots ,n$
   \begin{equation}\label{Ptheta.app}
    \ds x_s(t)=\sum_{ \tiny \begin{array}{c}
        \chi \in I_\theta  \end{array} }
    \; c^\chi_{s} \eu^{-(\chi_1 \wedge_1+ \ldots + \chi_n \wedge_n)t}+ o( \eu^{-(\Theta+\ep) t}) \,\,
    \textrm{ with $\chi=(\chi_1, \ldots , \chi_n)$}
  \end{equation}
  Moreover, let us denote by $d_i=c^\chi$ where $\chi$ is such that $\chi_{i}=1$ and $\chi_j=0$ for $j \ne i$. Then the coefficients
  $d_{i}$ just depend on $x(0)$ and the coefficients $c^\chi$
 where $|\chi|:=\sum \chi_i>1$ are determined by the coefficients $d_i$. In fact, for any $\chi$ such that $\chi_j=0$
   for $j \ge i+1$, then $c^\chi$ is determined by $d_1, \ldots , d_i$.
 \end{prop}
 Lemma \ref{nonres2} is a straightforward consequence of Lemma \ref{nonres}. The property of the coefficients
 is a straightforward consequence of the construction.

 From the discussion just after the proof of Lemma \ref{l:step.a} we get the following result,
  which proves Lemma \ref{datogliere} (in the non-resonant case).
  \begin{remark}\label{spiego.datogliere}
  The coefficients $d_i$ (from which all the other coefficients can be determined exactly) can be written as
  $d_i=P_i (x(0))- \tilde{\mathfrak{N}}_i $. However notice that they are the real unknown of this procedure, since
  $\mathfrak{N}_i$ is defined by a fixed point argument and we can just say that $d_i=P_i (x(0)) +O(x(0)^2)$.
 \end{remark}
 \begin{remark}\label{regolarita}
   In fact in the previous Lemma we can allow $N$ to be just $C^{\mu}$ where
   $\mu=\max\{|\chi|= \sum_i \chi_i \mid \chi=(\chi_1, \ldots , \chi_n) \in I_{\Theta} \}$
 \end{remark}

  If we drop the assumptions $\bs{R1}$, $\bs{R2}$ we have to deal with resonances, which at the end forces us to replace constants
  $c^\chi_s$ by polynomials.
  Assume first  $\bs{R1}$ but drop $\bs{R2}$. It might happen, e.g., that $\wedge_3=4 \wedge_2$: In this case in the functions
   $\tilde{M}_i^j(t)$ defined   in \eqref{redefinebis} we find terms of type $c \eu^{-4\wedge_2t}=c \eu^{- \wedge_3t}$,
   which are in resonance with $\eu^{L t}$. Thus
  in   $\tilde{a}_i^j(t)$   in \eqref{redefinebis} we find resonant terms of the form  $c t\eu^{- \wedge_3t}$.

  Further, if $\bs{R1}$ does not hold, the functions  $\tilde{a}^0_i(t)$ satisfies just the estimate in Remark \ref{rem1}.
  We assume that the eigenvalues $-\wedge_i$ are all real, this is the main case of interest in this article.
  It follows that each coordinate of  $\tilde{a}^0_i(t)$  equals  $q(t) \eu^{-\wedge_i t}$ where $q(t)$ is a polynomial of degree at
  most $k_i-1$ (and will depend on the coordinate). Then, plugging in the iterative scheme, we see that each coordinate
 $\tilde{a}^j_{i,s}(t)$ of
  $\tilde{a}^j_{i}(t)$ (for $s=1, \ldots, n$)  will be sum of exponentials, possibly multiplied by polynomials.
With this observation in mind  we get the following adaption of Lemma \ref{nonres2}.
    \begin{prop}\label{resonant}
 Assume $N \in C^{\infty}$, and that all the eigenvalues of $L$ are real and negative. Let $x(t) \to 0$ as $t \to +\infty$; then for any $\Theta>0$
 and any $s=1, \ldots, n$ we can find polynomials
  $c^\chi_s(t)$, $\chi \in I_{\Theta}$ such that
   \begin{equation}\label{Ptheta.app2}
    \ds x_s(t)=\sum_{ \tiny \begin{array}{c}
        \chi\in I_\theta  \end{array} }
    \; c^\chi_s(t) \eu^{-(\chi_1 \wedge_1+ \ldots + \chi_m \wedge_m)t}+ o( \eu^{-\Theta t}) \,\,
    \textrm{ with $\chi=(\chi_1, \ldots , \chi_m)$}
  \end{equation}
 Again  denote by $d_i(t)=(c^{\chi}_{1}(t), \ldots, c^{\chi}_{n})$ where $\chi$ is
 such that $\chi_i=1$ and $\chi_j=0$ for $j \ne i$. Then the vector of polynomials
  $d_i(t)$ just depend on $x(0)$ and   the vector of polynomials $c^\chi(t)$
 where $|\chi|:=\sum \chi_i>1$ are determined by the  $d_i(t)$. Further if $\chi_j=0$
   for $j \ge i+1$ then $c^\chi(t)$ is determined by $d_1, \ldots , d_i$.
 \end{prop}
\begin{remark}\label{resonant.rem}
If $\bs{R2}$ does not hold then there are $i$ and $\chi=(\chi_1, \ldots, \chi_{m})$ such that $\wedge_i=\sum_{j=1}^{i-1} \chi_j \wedge_j $.
In this case the vectorial function $c^\chi(t)$ is in fact   constant for any $\chi$ such that $\chi_1 \wedge_1+ \ldots + \chi_{i-1} \wedge_{i-1}<\wedge_i$,
i.e. for any $\chi \in I_{\wedge_i-\ep}$.
\end{remark}
\begin{remark}\label{spiego.res}
We emphasize that, even if we drop $\bs{R1}$ and $\bs{R2}$, Remarks \ref{spiego.datogliere} and \ref{regolarita} continue to hold, by construction.
However $\eu^{Lt}$ may contain terms which are polynomials multiplied by exponentials.
\end{remark}
\begin{remark}\label{resonant.rem2}
A result analogous to Lemma \ref{resonant} can be obtained also in the case where the eigenvalues of $L$
may be complex and conjugate, say $\la_i=-\wedge_i \pm i \omega$ (with negative real parts, $\wedge_i>0$). In this case the corresponding functions $\tilde{a}_i^0(t)$
will be sum of terms of the form $\eu^{-\wedge_i t}[c' \cos(\omega t) + c'' \sin(\omega t)]$. Then
each coordinate of the function
$c^\chi(t)$ might be a polynomial possibly multiplied by exponential and  sinusoidal functions.
\end{remark}

We emphasize that all this discussion in this article is in fact applied to the case of a $2$-dimensional non-autonomous system \eqref{si.na}
which can be regarded as a $3$-dimensional autonomous system \eqref{si.naas}, which has
 the critical point $\bs{P^+}$. The eigenvalues of the linearization of \eqref{si.naas} in $\bs{P^+}$
 are $\la_1$, $\la_2$ and $-\gamma$, and are real and negative.

However the dimension corresponding to $\gamma$  is special, since the corresponding equation can be solved independently
from the others (the solution is in fact $\zeta(s)=c\eu^{-\gamma s}$), and does not depend on the initial conditions.
Let us go back to the notation of Section 2.4.
\begin{remark}\label{suPsi1}
Consider  a solution $\bs{y}(s,\tau; \Q, l_s)$ of \eqref{si.na} converging to $\bs{P}$, and
  assume $-\gamma \le \la_1$; Then the function $\Psi(s)$ defined in \eqref{defPsi} and Proposition \ref{sintetizzo} does not depend on
  $\Q$ (in fact it depends on $\tau$).
\end{remark}
\begin{remark}\label{suPsi2}
Consider  a solution $\bs{y}(s,\tau; \Q, l_s)$ of \eqref{si.na} converging to $\bs{P}$, and
  assume $-\gamma = \la_1$. If $\la_1<\la_2$ the coefficient $a$ in \eqref{expandy} depend on $\Q$
  while $\Psi(s)$ contains a resonant term of the form $c_{(0,0,1)} s \eu^{\la_1 s}$ which does not depend on $\Q$.
  If $\la_1=\la_2$ the coefficient $a$ in \eqref{expandybis} depends on $\Q$
  while $\Psi(s)$ contains a doubly resonant term of the form $c_{(0,0,1)} s^2 \eu^{\la_1 s}$ which does not depend on $\Q$.
\end{remark}
\begin{proof}
  Remark \ref{suPsi2} immediately  follows from the construction in this Appendix, see in particular Remark \ref{suPsi1}. To prove Remark \ref{suPsi2}
  recall that the determination of $c_{(0,0,1)}$ and $a$ in both the cases goes back to Lemma \ref{linn}, so it just
  involves the resolution of the linearization of \eqref{si.naas} in $\bs{\mathcal{P}}$.
  This latter problem at the end can be regarded as
  a linear differential equation with constant coefficients and a forcing term which is a resonant exponential.
  Then it trivially follows that the most resonant part just depends on the forcing term, so we get  the conclusion.
\end{proof}



\end{document}